\documentclass[a4paper,12pt,fleqn]{article}
\usepackage{amsmath,amsthm,amssymb,bm,cite,wrapfig,subfigure,multicol,multirow,
takumatCom,boites,framed,pifont,enumerate,picture,ifthen,fourier,longtable}
\usepackage[dvipdfm]{graphicx}
\usepackage[polutonikogreek, english]{babel}
\newcommand{\keyword}[1]{{{\textit{#1}}}}

\newtheorem{theo}{Theorem}[section]
\newtheorem{thm}[theo]{Theorem}

\newtheorem{defi}[theo]{Definition}

\theoremstyle{definition}
\newtheorem{rema}[theo]{Remark}
\newtheorem{exmp}[theo]{Example}
\theoremstyle{theorem}
\newtheorem{prob}[theo]{Problem}
\newtheorem{algo}[theo]{Algorithm}

\makeatletter

\@addtoreset{equation}{section}

\@addtoreset{figure}{section}

\@addtoreset{table}{section}
\makeatother

\newcommand{\subsubsubsection}{\@startsection{paragraph}{4}{\z@}%
  {1.5\Cvs \@plus.5\Cdp \@minus.2\Cdp}%
  {.5\Cvs \@plus.3\Cdp}%
  {\reset@font\normalsize\sffamily}
}

\begin{document}
\thispagestyle{empty}

\newpage

\setlength{\baselineskip}{0.6cm}

\thispagestyle{empty}

\begin{center}
	\Large{An Approximate Approach to E-optimal Designs 
	for Weighted Polynomial Regression 
	by Using Tchebycheff Systems and Orthogonal Polynomials}\\
	\mbox{}\\
	\mbox{}\\
	\large{Takuma Takeuchi, Hiroto Sekido}\\
	\mbox{}\\
	\large{Graduate School of Informatics, Kyoto University}
\end{center}

\begin{center}

%
%
%

\vspace*{1.0cm}
\textbf{Abstract}
\end{center}

In statistics, experimental designs are methods for making efficient experiments.
E-optimal designs are the multisets of experimental conditions which minimize
the maximum axis of the confidence ellipsoid of estimators.
The aim of this thesis is to propose a new algorithm for constructing E-optimal designs approximately
for weighted polynomial regression with a nonnegative weight function.

First, an algorithm to calculate E-optimal designs for weighted polynomial regression 
of particular weight functions is discussed.
Next a new algorithm for constructing E-optimal designs approximately is proposed. 
Notions of the Tchebycheff systems and orthogonal polynomials are used in the proposed algorithm.
Finally in this thesis, the results of numerical examples are shown 
in order to verify the accuracy of the E-optimal designs computed by the proposed algorithm.

\newpage
\pagestyle{empty}
\tableofcontents
\setcounter{page}{0}
\newpage
\pagestyle{plain}


\section{Introduction}

In statistics, experimental designs are methods for making efficient experiments. 
Experimental designs are needed especially for experimenters. 
At first, in 1920s, Fisher \cite{Fisher1925} considered a formal mathematical methodology for designing experiments. 
This is the beginning of experimental designs.
We can make efficient experiments by analyzing a relationship between experimental conditions and the accuracy of estimators. 
In experimental designs, optimal designs are multisets of experimental conditions 
which give us the highest accuracy estimators based on a particular optimality criterion. 
On different demands of experimenters, several optimal criteria and optimal designs 
\cite{AtkinsonBogackaZhigljavsky2001, DetteStudden1997, Melas2005, Sekido2012}
are used.
One of the optimal criteria is the E-optimality criterion, which was introduced by Ehrenfeld \cite{Ehrenfeld1955}.
The best designs according to the E-optimality criterion are called E-optimal designs.
E-optimal designs minimize the maximum axis of the confidence ellipsoid of estimators, namely, 
E-optimal designs minimize the maximum eigenvalue of the covariance matrix of estimators. 
E-optimal designs have been investigated by numerous authors in the literature \cite{AtkinsonBogackaZhigljavsky2001, DetteStudden1997, Dette1993, Melas2005}.
E-optimal designs for only particular regression have been obtained exactly.

In this thesis, we discuss how to calculate E-optimal designs for weighted polynomial regression.
Weighted polynomial regression means polynomial regression with non-constant variance \cite{Fedorov1972}.
Optimal designs for weighted polynomial regression have been studied by many authors \cite{Studden1982, Dette1993, DetteHainesImhof1999}.

One of the approaches for obtaining E-optimal designs is to use the Tchebycheff systems.
One of the characteristics of the Tchebycheff systems is that 
there is a linear combination of their basis functions which satisfies some properties.
This linear combination is called the Tchebycheff function in this thesis.
The Tchebycheff systems play an important role in several domains of mathematics \cite{Karlin1968, KarlinStudden1966}.
For example, it is used for the theory of approximations, methods of interpolation, 
generalized moment problems, numerical analysis, oscillation properties of eigenfunctions of the Sturm--Liouville problems, 
generalized convexity, the theory of inequalities, and optimal designs.
E-optimal designs for polynomial regression and particular weighted polynomial regression
were studied with the Tchebycheff systems \cite{PukelsheimStudden1993, Dette1993}.
If Tchebycheff functions of the basis functions of the linear regression are known, 
then the corresponding E-optimal designs for general weighted polynomial regression can be calculated \cite{Melas2005}.
However, it is not trivial how to obtain Tchebycheff functions from the Tchebycheff systems.

Orthogonal polynomial sequences are families of polynomials such that 
the inner products of any distinct two polynomials in the sequences are zero.
Orthogonal polynomial sequences are also useful tools in many fields of mathematics \cite{Chihara1978, Szego1939, Nakamura2006}.
For example, it is used for the theory of approximations, and mathematical physics
including integrable systems.

In this thesis, we propose a new algorithm for constructing E-optimal designs approximately for 
weighted polynomial regression with general nonnegative 
weight functions
by using the Tchebycheff systems and orthogonal polynomials.
Moreover, we verify the accuracy of this algorithm by numerical examples.

Section \ref{sec: pre_ed} contains some preliminaries of optimal designs and the Tchebycheff systems.
Section \ref{sec: pre_op} contains some preliminaries of orthogonal polynomials and the Gram--Schmidt orthogonalization.
In Section \ref{sec: appEopt} we present an algorithm for constructing E-optimal designs approximately for 
weighted polynomial regression with general weight functions.
Section \ref{sec: numerical example} describes the results of numerical examples.
Section \ref{sec: concluding} is devoted to conclusions.

\section{Preliminaries of Experimental Designs} \label{sec: pre_ed}

\subsection{Linear Regression and Estimators}

A \keyword{linear regression model} is defined by
\begin{align}
	Y &= \theta^\top f(x) + \epsilon \notag \\
	&= 
	\begin{pmatrix}
		\theta_0, \theta_1, \cdots , \theta_{m - 1}
	\end{pmatrix}
	\begin{pmatrix}
		f_0(x) \\ f_1(x) \\ \vdots \\ f_{m - 1}(x)
	\end{pmatrix}
	+ \epsilon \label{eq: regression model}
\end{align}
where $f(x) = \left( f_0(x), f_1(x), \dots , f_{m - 1}(x) \right)^\top$ 
is a known vector of real-valued linearly independent continuous functions, 
$\theta = \left( \theta_0, \theta_1, \dots , \theta_{m - 1} \right)^\top$ 
is an unknown parameter vector, and $\epsilon$ is a random error term.
The functions $f_0(x), f_1(x), \dots ,$
$f_{m - 1}(x)$ are called \keyword{basis functions}.

The linear regression model \eqref{eq: regression model} means that 
$Y$ is the response of an observation at an experimental condition $x$.
The purpose in an experiment is to estimate the parameter vector $\theta$.
Here, we assume that all possible points where observations can be made are 
on the closed finite interval $\mathcal{X} = [a, b] \subset \rnum$.

Let us assume that we can make $N$ observations
\begin{align*}
	y_i = \theta^\top f(x_i) + \epsilon_i , \quad i = 1, 2, \dots , N 
\end{align*}
under the experimental conditions $x_1, x_2, \dots , x_N \in \mathcal{X}$. 
Throughout this thesis, we assume that the expectation of an error $\epsilon_i$ is zero and 
different errors are uncorrelated.
Conventionally we sometimes assume that the variance of an error is a positive constant.
That is, 
\begin{align}
	\ev [\epsilon_i] = 0, \quad \ev [\epsilon_i \epsilon_{j}] = 0, \quad
	\va [\epsilon_i] = \sigma^2 > 0, \quad i, j = 1, 2, \dots , N, \quad i \neq j. \label{eq: error assumptions}
\end{align}
The \keyword{best linear unbiased estimator} (\keyword{BLUE}) $\hat\theta$ of the parameter vector $\theta$
is defined as the estimator which satisfies the following three conditions:
\begin{enumerate}
	\item \label{enum: BLUE condition 1} The estimator $\hat\theta$ is described as a linear combination of the responses, 
	namely $\hat\theta = L \bm{y}$, where $L$ is an $m \times N$ matrix.
	\item \label{enum: BLUE condition 2} The expectation of the estimator $\hat\theta$ is equal to $\theta$, 
	namely $\ev \left[ \hat\theta \right] = \theta$.
	\item For an arbitrary estimator $\bar\theta$ which satisfies the conditions 
	\eqref{enum: BLUE condition 1} and \eqref{enum: BLUE condition 2}, 
	$\cov \left[ \hat\theta \right] - \cov \left[ \bar\theta \right]$ is nonnegative definite, 
	where $\cov \left[ \theta \right]$ denotes the covariance matrix of $\theta$.
\end{enumerate}

The following theorem is well known in statistics.
\begin{thm}[Gauss--Markov's Theorem] \label{thm: Gauss--Markov theorem}
	Under the conditions \eqref{eq: error assumptions} and \\ $\det (X^\top X) \neq 0$, 
	the BLUE $\hat\theta$ of the parameter vector $\theta$ is given by
	\begin{align}
		\hat\theta = \left( X^\top X \right)^{-1} X^\top y , \label{eq: gaussmarkov BLUE}
	\end{align}
	where $X = \left( f(x_1), f(x_2), \dots , f(x_N) \right)^\top$ is an $N \times m$ matrix, 
	and $y = (y_1, y_2, \dots , y_N)^\top$.
	The covariance matrix of the BLUE $\hat\theta$ is given by
	\begin{align*}
		\cov \left[ \hat\theta \right] = \sigma^2 \left( X^\top X \right)^{-1} .
	\end{align*} 
\end{thm}

\subsection{Optimal Designs and Fisher Information Matrix}

A \keyword{design} $\tilde \mu$ is a multiset of experimental conditions $x_1, x_2, \dots , x_N \in \mathcal{X}$.
When we make experiments, we should choose the optimal multiset ${\tilde \mu}^\ast$.
But in general, it is difficult to calculate the optimal multiset ${\tilde \mu}^\ast$.

Then, in this thesis we consider a multiset $\tilde \mu$ as a probability measure $\mu$.
Let $\mathcal{P}_\mathcal{X}$ denote the set of all probability measures on the Borel sets of the interval $\mathcal{X}$.
For given $\mu \in \mathcal{P}_\mathcal{X}$, let $\mu (x)$ denote the cumulative distribution function, and let 
$\mathrm{Prob}_\mu (x)$ denote the probability mass function 
\begin{align*}
	\mathrm{Prob}_\mu (x) = \lim_{\delta \to +0} \left( \mu (x + \delta ) - \mu (x - \delta ) \right) .
\end{align*}
We consider the probability measure whose the probability mass function is given by 
\begin{align*}
	\mathrm{Prob}_\mu (x) = \frac{\# \left\{ i \in \left\{ 1, 2, \dots , N \right\} \big| x_i = x \right\}}{N} , 
\end{align*}
where $\# S$ denotes the number of elements in the set $S$.
Assume that the distinct points among $x_1, x_2, \dots , x_N$ are 
the points $x_1, x_2, \dots , x_n, \; n \leq N$.
Thus we also call the probability measure $\mu$ the \keyword{design}, namely, 
the design $\mu$ means that we make $N \rho_i$ experiments under a condition $x_i$, \; $i = 1, 2, \dots , n$, 
where 
\begin{align*}
	\rho_i = \mathrm{Prob}_\mu (x_i) , \quad i = 1, 2, \dots , n .
\end{align*}
We sometimes write the design $\mu$ as
\begin{align}
	\mu = 
	\begin{pmatrix}
		x_1 & x_2 & \cdots & x_n \\
		\rho_1 & \rho_2 & \cdots & \rho_n
	\end{pmatrix} . \label{eq: design having a finite number of supports}
\end{align}

For a fixed sample size $N$, 
let us consider the case where 
the numbers $N \rho_i,$ $i = 1, 2, \dots , n$ are not necessary to be integers.
That is, $\rho_i, \; i = 1, 2, \dots , n$ must be arbitrary nonnegative numbers such that $\sum_{i = 1}^n \rho_i = 1$.
In practice, the numbers $N \rho_i, \; i = 1, 2, \dots , n$ of the design $\mu$ are sometimes rounded to be integers
in order to consider the corresponding multiset $\tilde \mu$.
Thus, hereinafter the design $\mu$ denotes only a probability measure, not a multiset.

\begin{figure}[!h]
	\begin{center}
		\includegraphics[width=60mm, clip]{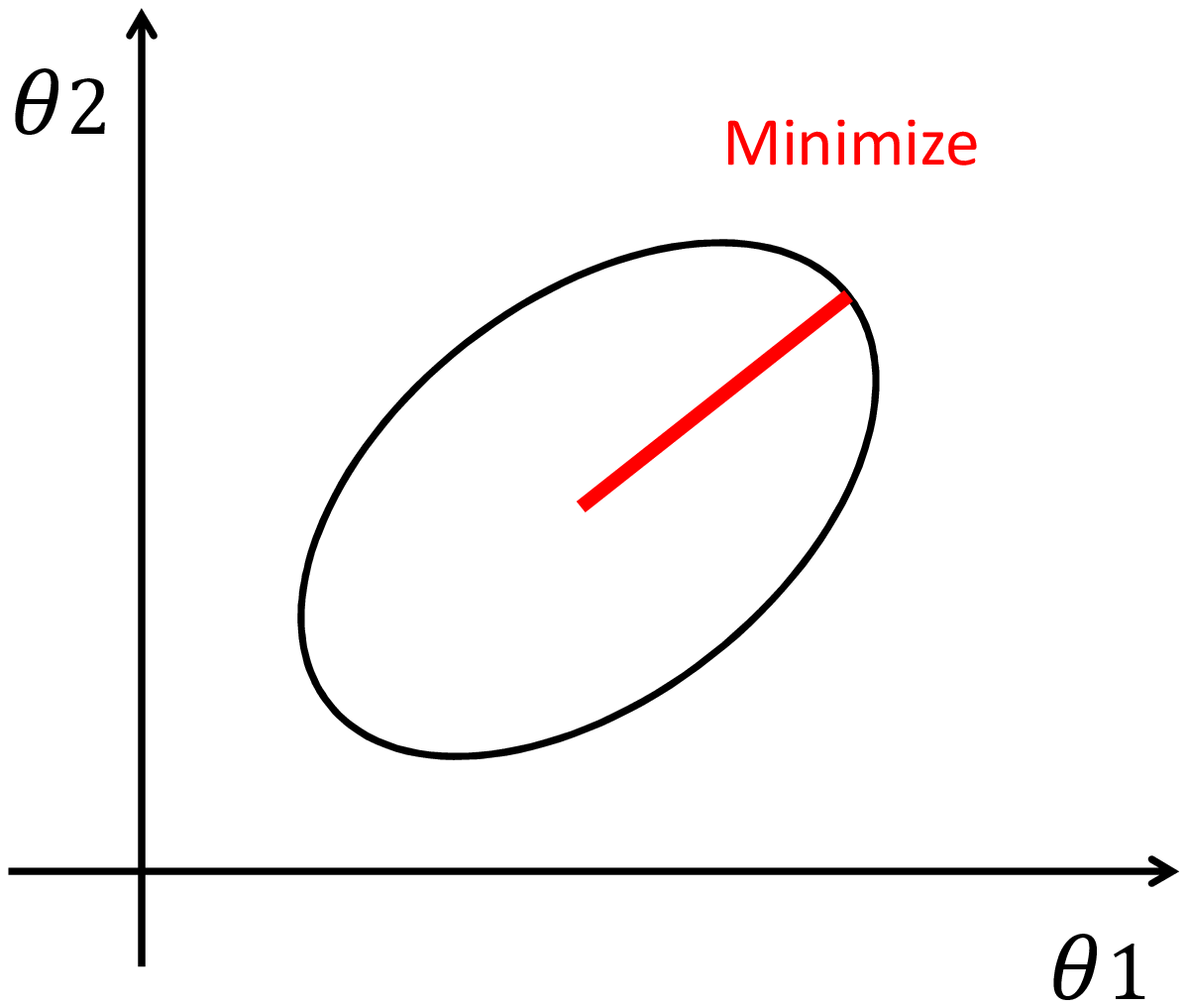} 
		\caption{The characteristics of E-optimal designs} \label{fig: E-optimal design}
	\end{center}
	\begin{center}
		\includegraphics[width=60mm, clip]{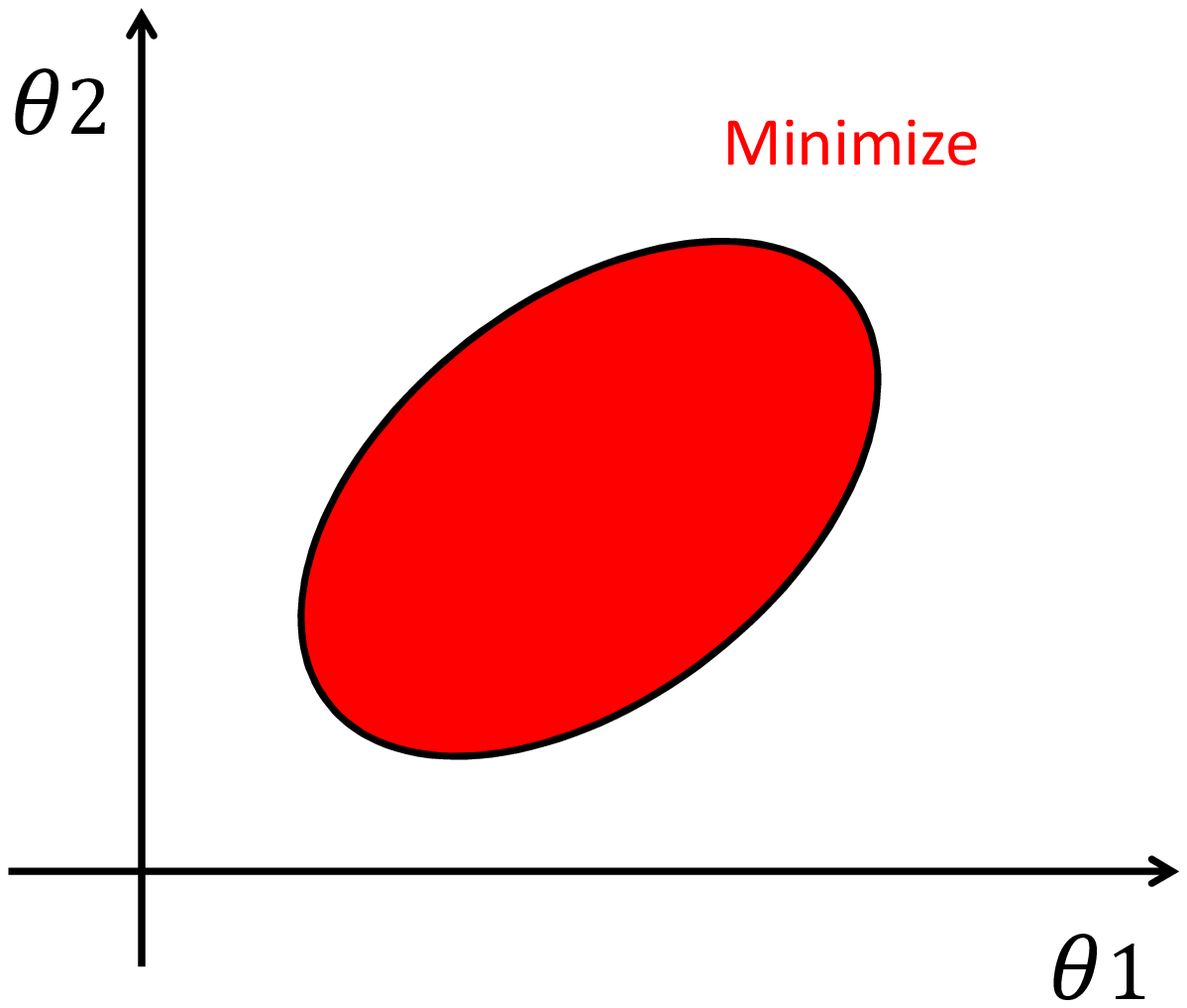} 
		\caption{The characteristics of D-optimal designs} \label{fig: D-optimal design}
	\end{center}
	\begin{center}
		\includegraphics[width=60mm, clip]{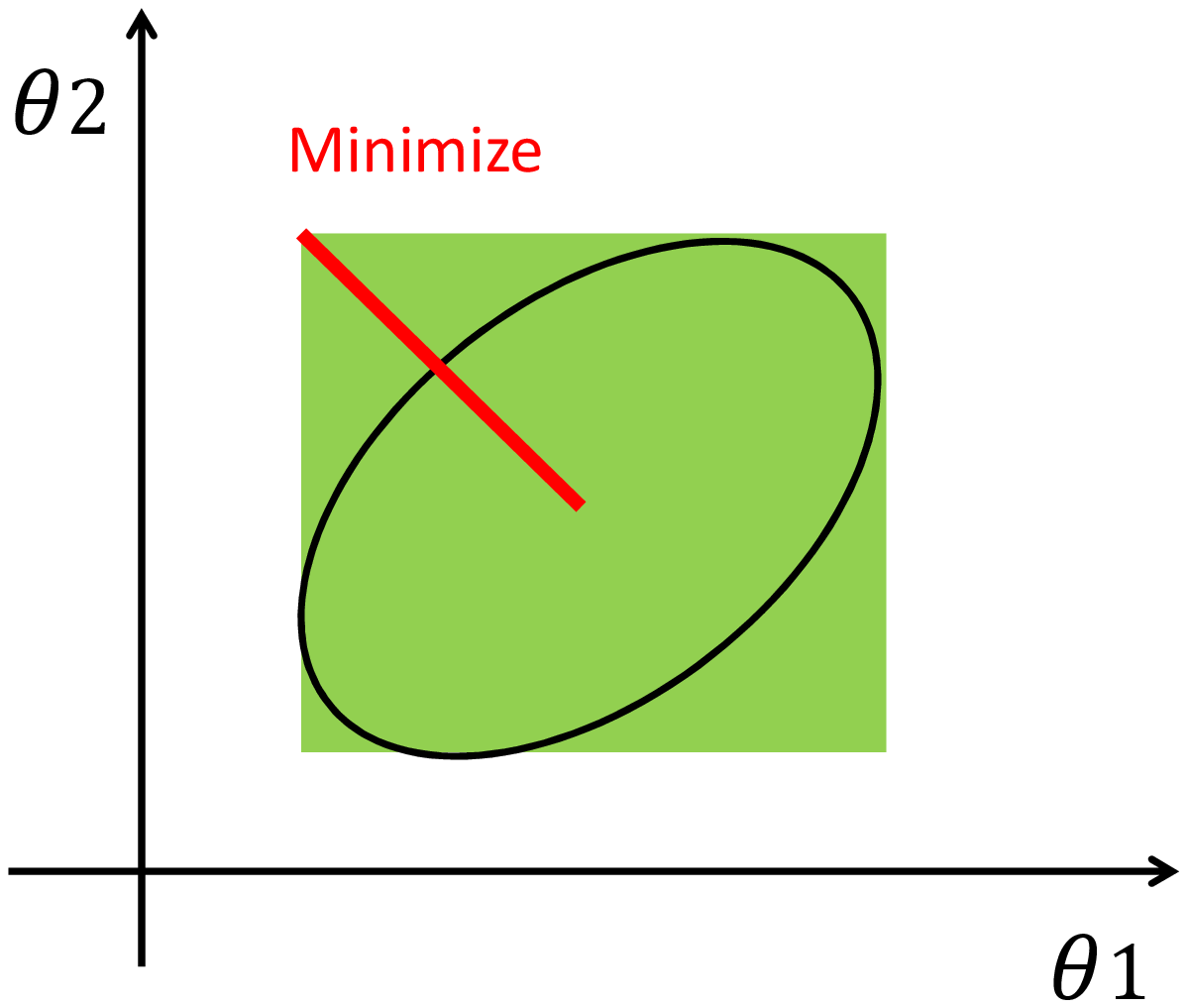} 
		\caption{The characteristics of A-optimal designs} \label{fig: A-optimal design}
	\end{center}
\end{figure}

We should choose a good design, 
since the BLUE $\hat\theta$ \eqref{eq: gaussmarkov BLUE} depends on the design $\mu$.
In general, if the covariance matrix of the BLUE $\hat\theta$ is ``{\textit{small}}'' in some sense, 
the BLUE $\hat\theta$ becomes a highly accurate estimator. 
Here, in order to define what means that the covariance matrix is small, 
let us consider the \keyword{Fisher information matrix}.
The Fisher information matrix of the design $\mu$ is defined by the Gram matrix
\begin{align}
	M(\mu ) &= \int_\mathcal{X} f(x) f^\top (x) \di \mu (x) \label{eq: def of fisher information matrix} \\
	&= \sum_{i = 1}^n f(x_i) f^\top (x_i) \rho_i \notag \\
	&=
	\begin{pmatrix}
		\sum_{i = 1}^n f_0 (x_i) f_0 (x_i) \rho_i & \dots & \sum_{i = 1}^n f_0 (x_i) f_{m - 1} (x_i) \rho_i \\
		\vdots & \ddots & \vdots \\
		\sum_{i = 1}^n f_{m - 1} (x_i) f_0 (x_i) \rho_i & \dots & \sum_{i = 1}^n f_{m - 1} (x_i) f_{m - 1} (x_i) \rho_i \\
	\end{pmatrix} \notag .
\end{align}
By Theorem \ref{thm: Gauss--Markov theorem}, the covariance matrix of the BLUE $\hat\theta$ is represented as
\begin{align}
	\cov \left[ \hat\theta \right] = \frac{\sigma^2}{N} M^{-1}(\mu ) . \label{eq: Gauss--Markov result}
\end{align}
In order to make the covariance matrix of the BLUE $\hat\theta$ the {\textit{smallest}} in some sense,
we should choose the \keyword{optimal design} $\mu$ 
whose Fisher information matrix $M(\mu )$ takes the ``{\textit{smallest form}}''.
Here, let us consider the \keyword{$\Phi_p$-optimality criterion}, 
a commonly used optimality criterion in experimental designs, that is
\begin{align}
	\begin{split}
		& \underset{\mu}{\text{Minimize}} \ \Phi_p (\mu ) = \left( \frac{1}{m} \tr M^{-p} (\mu ) \right)^{\frac{1}{p}}, \quad 0 < p < \infty \\
		& \text{subject to} \ \mu \in \mathcal{P}_\mathcal{X} . 
	\end{split}
	\label{eq: def phi_p}
\end{align}
Especially, when $p \to \infty, \ p \to 0, \ p = 1$, the objective function $\Phi_p (\mu )$ is represented as
\begin{align*}
	& \Phi_{\infty} (\mu ) = \underset{1 \leq i \leq m}{\max} \frac{1}{\lamb_i (\mu )} , \\
	& \Phi_{0} (\mu ) = \left( \det M(\mu ) \right)^{\frac{1}{m}} , \\
	& \Phi_{1} (\mu ) = \frac{1}{m} \tr M^{-1}(\mu ) , 
\end{align*}
respectively, where $\lamb_1 (\mu ), \lamb_2 (\mu ), \dots , \lamb_m (\mu )$ denote 
the eigenvalues of the Fisher information matrix $M(\mu )$ of the design $\mu$.
The design $\mu$ which minimizes $\Phi_{\infty} (\mu )$ is called an \keyword{E-optimal design}.
Similarly, the design $\mu$ which minimizes $\Phi_{0} (\mu )$ is called a \keyword{D-optimal design},
and the design $\mu$ which minimizes $\Phi_{1} (\mu )$ is called an \keyword{A-optimal design}.
That is, the E-optimal designs are the optimal solutions of the optimization problem
\begin{align}
	\begin{split}
		& \underset{\mu}{\text{Maximize}} \ \lamb_{\min} (M (\mu )) \\ 
		& \text{subject to } \ \mu \in \mathcal{P}_{\mathcal{X}} , 
	\end{split}
	\label{eq: the optimization problem of E-optimal designs}
\end{align}
where $\lamb_{\min} (M (\mu ))$ denotes the minimum eigenvalue of the Fisher information matrix $M (\mu )$.
The D-optimal designs are the optimal solutions of the optimization problem
\begin{align}
	\begin{split}
		& \underset{\mu}{\text{Maximize}} \ \det M(\mu ) \\
		& \text{subject to } \ \mu \in \mathcal{P}_{\mathcal{X}} , 
	\end{split}
	\label{eq: the optimization problem of D-optimal designs}
\end{align}
and the A-optimal designs are the optimal solutions of the optimization problem
\begin{align}
	\begin{split}
		& \underset{\mu}{\text{Minimize}} \ \tr M^{-1}(\mu ) \\
		& \text{subject to } \ \mu \in \mathcal{P}_{\mathcal{X}} . 
	\end{split}
	\label{eq: the optimization problem of A-optimal designs}
\end{align}
When the error $\epsilon$ is normally distributed, 
the \keyword{confidence ellipsoid} for the BLUE $\hat\theta$ with an arbitrary fixed confidence level
is defined by
\begin{align}
	\left\{ \bar\theta \ \big | \ \left( \bar\theta - \hat\theta \right)^\top M^{-1}(\mu ) \left( \bar\theta - \hat\theta \right) \leq c \right\} , \label{eq: def of confidence ellipsoid}
\end{align}
where $c$ is a constant depending only on the confidence level.
In this case E-, D-, and A-optimal designs can be interpreted geometrically in terms of the confidence ellipsoid.
E-optimal designs minimize the size of the major axis of the confidence ellipsoid, 
D-optimal designs minimize the volume of the confidence ellipsoid, 
and A-optimal designs minimize the dimension of the diagonal of the enclosing box around the confidence ellipsoid respectively. 
These characteristics are shown by 
Figures \ref{fig: E-optimal design}, \ref{fig: D-optimal design}, and \ref{fig: A-optimal design} .

\clearpage

In this thesis, we discuss the problem of calculating E-optimal designs for the linear regression on $\mathcal{X} = [-1, 1]$.
\begin{breakbox}
	\begin{prob}[The problem of calculating E-optimal designs] \label{prob: E-optimal design}
		We consider the linear regression \eqref{eq: regression model}.
		The E-optimal designs are the optimal solutions of the optimization problem
		\begin{align*}
			& \underset{\mu}{\text{Maximize}} \ \lamb_{\min} (M (\mu )) \\ 
			& \text{subject to } \ \mu \in \mathcal{P}_{[-1, 1]}, 
		\end{align*}
		where $\mu$ is a design, $M(\mu )$ is a Fisher information matrix for the regression model
		defined by \eqref{eq: def of fisher information matrix}, and 
		$\mathcal{P}_{[-1, 1]}$ is the set of all probability measures on $[-1, 1]$.
	\end{prob}
\end{breakbox}

We note that E-optimal designs do not depend on the sample size $N$.

\subsection{Tchebycheff Systems and Their Applications to Optimal Designs}

Let $u_1, u_2, \dots , u_n: I \to \rnum$ denote linearly independent continuous functions
defined on a closed finite interval $I = [a, b]$.
If the determinant
\begin{align}
	\begin{vmatrix}
		u_1 (t_1) & u_1 (t_2) & \cdots & u_1 (t_n) \\
		u_2 (t_1) & u_2 (t_2) & \cdots & u_2 (t_n) \\
		\vdots & \vdots & \ddots & \vdots \\
		u_n (t_1) & u_n (t_2) & \cdots & u_n (t_n) 
	\end{vmatrix}
	\label{eq: determinant of Tsys}
\end{align}
is always positive or always negative whenever the parameters are chosen as $a \leq t_1 < t_2 < \dots < t_n \leq b$, 
the set $\left\{ u_1, u_2, \dots , u_n \right\}$ is called the \keyword{Tchebycheff system} on $I$.
If the determinant \eqref{eq: determinant of Tsys} is always nonnegative or always nonpositive 
whenever the parameters are chosen as $a \leq t_1 < t_2 < \dots < t_n \leq b$, 
then the set $\left\{ u_1, u_2, \dots , u_n \right\}$ of the functions is called the \keyword{weak Tchebycheff system} on $I$.
It is well known \cite[Theorem II 10.2]{KarlinStudden1966}
that if the set $\{ u_1, u_2, \dots , u_n \}$ is a weak Tchebycheff system, 
then there exists a unique function $\kappa (t)$ given by
\begin{align*}
	\kappa (t) = \gamma^\top u(t), \quad \gamma = \left( \gamma_1, \gamma_2, \dots , \gamma_n \right)^\top \in \rnum^n, 
	\quad u(t) = \left( u_1, u_2, \dots , u_n \right)^\top 
\end{align*}
which satisfies the following properties: 
\begin{enumerate}
	\item $\left| \kappa (t) \right| \leq 1$ for all $t \in I$, 
	\item There exist $n$ points $s_1, s_2, \dots , s_n$ chosen as $a \leq s_1 < s_2 < \dots < s_n \leq b$ such that
	$\kappa (s_i) = (-1)^i, \; i = 1, 2, \dots , n$. 
\end{enumerate}
The function $\kappa (t)$ is called the \keyword{Tchebycheff function}, 
and the points $s_1, s_2, \dots , s_n$ are called the \keyword{Tchebycheff points}.
\begin{exmp}[Examples of Tchebycheff systems] \label{exmp: examples of Tchebycheff systems}
	In each of the following cases, the set $\{ u_1, u_2, \dots , u_n \}$ is 
	a Tchebycheff system \cite[pp. 9--20]{KarlinStudden1966}.
	\begin{enumerate}
		\item \label{item: examples of Tchebycheff systems power function} 
		The set of power functions $u_i(t) = h(t) t^{i} , \; i = 1, 2, \dots , n$ is a Tchebycheff system 
		where $h(t)$ is a positive continuous function on $I$.
		If $h(t)$ is a nonnegative continuous function on $I$, 
		then the set $\{ u_1, u_2, \dots , u_n \}$ is a weak Tchebycheff system.
		\item The set of eigenfunctions of the Sturm--Liouville operator 
		\begin{align*}
			L(\phi ) = - \dif{}{t} \left( p(t) \dif{\phi}{t} \right) + q(t) \phi
		\end{align*}
		is a Tchebycheff system 
		where $p(t)$ is a continuous positive function on $I$, 
		and $q(t)$ is a continuous function on $I$.
	\end{enumerate}
\end{exmp}

There are some relations between the Tchebycheff systems and optimal designs.
The following theorem indicates the way to obtain E-optimal designs by using the Tchebycheff systems \cite[pp. 94-97]{Melas2005}.
\begin{thm}[Tchebycheff designs] \label{thm: Tchebycheff design}
	Suppose the set 
	$\left\{ f_0(x), f_1(x), \dots , f_{m - 1}(x) \right\}$ 
	of the basis functions of the linear regression \eqref{eq: regression model}
	is a Tchebycheff system and 
	generates a Tchebycheff function $\kappa (x)$ given by
	\begin{align*}
		\kappa (x) = \gamma^\top f(x), \quad \gamma \in \rnum^n, \quad x \in \mathcal{X}. 
	\end{align*}
	Also suppose that $s_1, s_2, \dots , s_n$ are Tchebycheff points of the Tchebycheff function $\kappa (x)$.
	If Tchebycheff points of the system are determined uniquely, 
	then the design $\mu^\ast$ given by
	\begin{align*}
		& \mu^{\ast} = 
		\begin{pmatrix}
			s_1 & s_2 & \dots & s_m \\
			\rho_1 & \rho_2 & \dots & \rho_m
		\end{pmatrix}
		, \\
		& \left( \rho_1, \rho_2, \dots , \rho_m \right)^\top = \frac{F^{-1} \gamma}{\gamma^\top \gamma} , 
		\quad F = \left( f_{i - 1} (s_j) \cdot (-1)^{j + 1} \right)_{i, j = 1, 2, \dots , m}
	\end{align*}
	is called the \keyword{Tchebycheff design}.
	The Tchebycheff design $\mu^\ast$ is the E-optimal design 
	if the linear regression \eqref{eq: regression model}
	has a unique E-optimal design.
\end{thm}

\subsection{Optimal Designs for Some Regression Models}

\subsubsection{Optimal Designs for Polynomial Regression}

Let us consider the linear regression model \eqref{eq: regression model}.
In the case where $f(x)$ is given by 
\begin{align*}
	f(x) = \left( f_0(x), f_1(x), \dots , f_{m - 1}(x) \right)^\top = \left( 1, x, \dots , x^{m - 1} \right)^\top ,
\end{align*}
the linear regression model 
\begin{align}
	y = \sum_{i = 0}^{m - 1} \theta_i x^i + \epsilon \label{eq: polynomial regression}
\end{align}
is called a \keyword{polynomial regression model}.
Let us take $N$ observations
\begin{align}
	y_k = \sum_{i = 0}^{m - 1} \theta_i x_k^i + \epsilon_k, \quad k = 1, 2, \dots , N \label{eq: polynomial regression}
\end{align}
under the experimental conditions $x_1, x_2, \dots , x_N \in \mathcal{X}$. 
Several researches show how to compute optimal designs for polynomial regression.
For example, D-optimal designs for polynomial regression can be calculated 
by using canonical moments \cite{DetteStudden1997}.
E-optimal designs for polynomial regression can be calculated 
by using the Tchebycheff systems \cite{Dette1993}.

\subsubsection{Optimal Designs for Weighted Polynomial Regression}

Let us consider the polynomial regression model \eqref{eq: polynomial regression} 
without the assumption that the variance of an error is constant.
Namely, we assume that
\begin{align}
	\ev [\epsilon_i] = 0, \quad \ev [\epsilon_i \epsilon_j] = 0 , \quad
	\va [\epsilon_i] = \frac{\sigma^2}{\omega (x_k)} > 0, \quad 
	i, j = 1, 2, \dots , N , \quad i \neq j, \label{eq: weighted error assumptions}
\end{align}
where $\omega (x)$ is a nonnegative function, called a \keyword{weight function of regression},
which depends on an experimental condition $x$.
The Fisher information matrix of a design $\mu$ for weighted polynomial regression 
is redefined as
\begin{align*}
	M(\mu ) &= \int_\mathcal{X} \omega(x) f(x) f^\top (x) \di \mu \\
	&= \int_\mathcal{X} \left( \sqrt{\omega(x)} f(x) \right) \left( \sqrt{\omega(x)} f(x) \right)^\top \di \mu .
\end{align*}
Then, E-, D-, and A-optimal designs are defined by the same ways 
\eqref{eq: the optimization problem of E-optimal designs}, 
\eqref{eq: the optimization problem of D-optimal designs}, and 
\eqref{eq: the optimization problem of A-optimal designs}
respectively, as the polynomial regression model.

E-optimal designs for weighted polynomial regression can be calculated exactly for only particular weight functions.
These are described later in Subsection \ref{sec: e-optimal design in dette weight function}.

%

%

\section{Preliminaries of Orthogonal Polynomials} \label{sec: pre_op}

In this section we discuss orthogonal polynomials and the Gram--Schmidt orthogonalization.
The Gram--Schmidt orthogonalization is a method for making an orthogonal polynomial sequence \cite{Chihara1978, Szego1939, Nakamura2006}.

\subsection{Inner Product, Moments, and Classical Orthogonal Polynomials}

At first, we define the inner product $\langle \cdot , \cdot \rangle$ 
with respect to a nonnegative function $\eta (x)$ by
\begin{align}
	\inpro{p(x)}{q(x)} = \int_a^b p(x) q(x) \eta (x) \di x \label{eq: def inner product} ,
\end{align}
where $p, q$ are polynomials defined on $[a, b]$.
The function $\eta (x)$ is called a \keyword{weight function of orthogonal polynomials}.

If polynomials $p(x), q(x)$ satisfy 
\begin{align*}
	\inpro{p(x)}{q(x)} = 0 , 
\end{align*}
then we say that they are \keyword{orthogonal}.
The following functions are examples of orthogonal polynomials. They are called 
the \keyword{classical orthogonal polynomials}.
\begin{enumerate}
	\item Jacobi polynomials $\jacobiP_n^{(\alpha , \beta )}(x)$:
	\begin{align}
		\jacobiP_n^{(\alpha , \beta )}(x) 
		&= \frac{(-1)^n}{2^n n! (1 - x)^\alpha (1 + x)^\beta} \dif{^n}{x^n} \left( (1 - x)^{n + \alpha} (1 + x)^{n + \beta} \right) \notag \\
		&= \frac{1}{2^n} \sum_{k = 0}^n \ve{n + \alpha}{n - k} \ve{n + \beta}{k} (x - 1)^k (x + 1)^{n - k}, \label{eq: Jacobi polynomial}
	\end{align}
	where $n \in \znum_{\geq 0}$, and $\alpha > -1, \beta > -1$.
	The orthogonality relation is given by
	\begin{align*}
		& \int_{-1}^1 \jacobiP_m^{(\alpha , \beta )}(x) \jacobiP_n^{(\alpha , \beta )}(x) \eta_{\mathrm{\jacobiP}}^{(\alpha , \beta )} (x) \di x \\
		& \quad = \frac{2^{\alpha + \beta + 1} \Gamma (n + \alpha + 1) \Gamma (n + \beta + 1)}{(2 n + \alpha + \beta + 1) \Gamma (n + \alpha + \beta + 1) n!} \delt_{m, n}
	\end{align*}
	where the weight function $\eta_{\mathrm{\jacobiP}}^{(\alpha , \beta )} (x)$ is given by $\eta_{\mathrm{\jacobiP}}^{(\alpha , \beta )} (x) = (1 - x)^\alpha (1 + x)^\beta$.
	\item Laguerre polynomials $\laguerreP_n^{(\alpha )}(x)$:
	\begin{align*}
		\laguerreP_n^{(\alpha )}(x)
		&= \frac{\ex^x}{n! x^\alpha} \dif{^n}{x^n} \left( x^{n + \alpha} \ex^{-x} \right) \\
		&= \sum_{k = 0}^n \ve{n + \alpha}{n - k} \frac{(-x)^k}{k!}, 
	\end{align*}
	where $n \in \znum_{\geq 0}$, and $\alpha > -1$.
	The orthogonality relation is given by
	\begin{align*}
		\int_{0}^\infty \laguerreP_m^{(\alpha )}(x) \laguerreP_n^{(\alpha )}(x) \eta_{\mathrm{\laguerreP}}^{(\alpha )} (x) \di x
		= \frac{\Gamma (n + \alpha + 1)}{n!} \delt_{m, n}
	\end{align*}
	where the weight function $\eta_{\mathrm{\laguerreP}}^{(\alpha )} (x)$ is given by $\eta_{\mathrm{\laguerreP}}^{(\alpha )}  (x) = x^\alpha \ex^{-x}$.
	\item Hermite polynomials $\hermiteP_n (x)$: 
	\begin{align*}
		\hermiteP_n (x) 
		&= (-1)^n \ex^{x^2} \dif{^n}{x^n} \ex^{-x^2} \\
		&= n! \sum_{k = 0}^{\lfloor n / 2 \rfloor} \frac{(-1)^k (2 x)^{n - 2 k}}{(n - 2 k)! k!}, 
	\end{align*}
	where $n \in \znum_{\geq 0}$ and $\lfloor t \rfloor$ denotes the largest integer not exceeding $t$.
	The orthogonality relation is given by
	\begin{align*}
		\int_{-\infty}^\infty \hermiteP_m (x) \hermiteP_n (x) \eta_{\mathrm{\hermiteP}}(x) \di x
		= 2^n n! \sqrt{\pie} \delt_{m, n}
	\end{align*}
	where the weight function $\eta_{\mathrm{\hermiteP}} (x)$ is given by $\eta_{\mathrm{\hermiteP}} (x) = \ex^{-x^2}$.
\end{enumerate}

\subsection{Gram--Schmidt Orthogonalization}

The following algorithm, called the Gram--Schmidt orthogonalization, means a method for orthogonalizing a set of polynomials
in an inner product space.
Here we consider the inner product space defined by \eqref{eq: def inner product}.

\begin{algo}[Gram--Schmidt orthogonalization] \label{algo: Gram--Schmidt}
	We define the projection operator $\mathrm{proj}$ by 
	\begin{align*}
		\mathrm{proj}_{v} \left( u \right)
		= \frac{\inpro{u}{v}}{\inpro{v}{v}} v,
	\end{align*}
	where $u, v$ are polynomials.
	Then, if $u_1, u_2, \dots , u_n$ are linearly independent polynomials, 
	the following process
	\begin{align}
		& \mathrm{for} \ k = 1, 2, \dots , n: \notag \\
		& \quad v_k = u_k - \sum_{l = 1}^{k - 1} \mathrm{proj}_{v_l} \left( u_k \right) \label{eq: Gram--Schmidt 1}
	\end{align}
	constructs orthogonal polynomials $v_1, v_2, \dots , v_n$.
\end{algo}

The result $v$ of the Gram--Schmidt orthogonalization can be expressed by
\begin{align}
	v_k = \frac{
	\begin{vmatrix}
		\inpro{u_1}{u_1} & \inpro{u_1}{u_2} & \dots  & \inpro{u_1}{u_k} \\
		\inpro{u_2}{u_1} & \inpro{u_2}{u_2} & \dots  & \inpro{u_2}{u_k} \\
		\vdots           & \vdots           & \ddots & \vdots           \\
		\inpro{u_{k - 1}}{u_1} & \inpro{u_{k - 1}}{u_2} & \dots  & \inpro{u_{k - 1}}{u_k} \\
		u_1              & u_2              & \dots  & u_k
	\end{vmatrix}
	}
	{
	\begin{vmatrix}
		\inpro{u_1}{u_1} & \inpro{u_1}{u_2} & \dots  & \inpro{u_1}{u_{k - 1}} \\
		\inpro{u_2}{u_1} & \inpro{u_2}{u_2} & \dots  & \inpro{u_2}{u_{k - 1}} \\
		\vdots           & \vdots           & \ddots & \vdots           \\
		\inpro{u_{k - 1}}{u_1} & \inpro{u_{k - 1}}{u_2} & \dots  & \inpro{u_{k - 1}}{u_{k - 1}} 
	\end{vmatrix}
	} 
	, \quad k = 1, 2, \dots , n. \label{eq: Gram--Schmidt determinant}
\end{align}

\section{Construction of E-optimal Designs for Weighted Polynomial Regression} \label{sec: appEopt}

Let us consider about the weighted polynomial regression with the weight function $w(x)$ of regression.
From Example \ref{exmp: examples of Tchebycheff systems} \eqref{item: examples of Tchebycheff systems power function}, 
the set 
\begin{align}
	\left\{ \sqrt{w(x)}, x \sqrt{w(x)}, \dots , x^{m - 1} \sqrt{w(x)} \right\} \label{eq: the set of weighted function}
\end{align}
is a Tchebycheff system.
By Theorem \ref{thm: Tchebycheff design}, 
if we know the Tchebycheff function of the set \eqref{eq: the set of weighted function}, 
we can obtain E-optimal designs for weighted polynomial regression.
However, it is not trivial how to get the Tchebycheff functions.
In this section, we propose an approximate approach to construct E-optimal designs for weighted polynomial regression.

\subsection{Construction of E-optimal Designs for Particular Weighted Polynomial Regression with Jacobi Polynomials} \label{sec: e-optimal design in dette weight function}

In this subsection, we discuss an example 
that we can compute E-optimal designs exactly for weighted polynomial regression 
using the Tchebycheff systems.
Let us consider weighted polynomial regression in the case where 
its weight function is described by
\begin{align}
	w(x) = (1 - x)^\alpha (1 + x)^\beta , \quad \alpha , \beta \in \left\{ 0, 1 \right\} . \label{eq: dette weight function}
\end{align}
Dette \cite{Dette1993} shows how to compute E-optimal designs for these regression as the following theorem.

\begin{thm} \label{thm: dette E-optimal design for special weight function}
	With the weight function \eqref{eq: dette weight function}, 
	the function
	\begin{align}
		(1 - x)^{\alpha / 2} (1 + x)^{\beta / 2} J_{m - 1}^{\left( \alpha - 1 / 2, \ \beta - 1 / 2 \right)} \label{eq: yamate}
	\end{align}
	is a Tchebycheff function where $J_m^{(\alpha , \beta )}$ denotes a Jacobi polynomial \eqref{eq: Jacobi polynomial}. 
	Then, the Tchebycheff design $\mu^\ast$ for weighted polynomial regression with the weight function $w(x)$
	can be computed by Theorem \ref{thm: Tchebycheff design}.
	The Tchebycheff design $\mu^\ast$ is equal to the E-optimal design.
\end{thm}

\begin{figure}[!h]
	\begin{center}
		\includegraphics[width=60mm, clip]{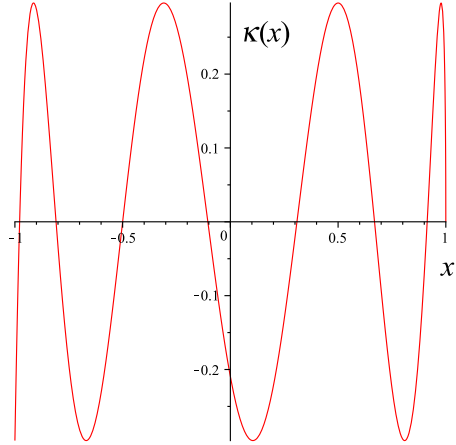} \\
		\caption{A Tchebycheff function $\kappa (x) = J_7^{(1/2, -1/2)} (x) \sqrt{1 - x}$} \label{fig: jacobijaco}
	\end{center}
\end{figure}

Figure \ref{fig: jacobijaco} shows the graph of \eqref{eq: yamate} with $\alpha = 0, \; \beta = 1, \; m = 8$.
This figure indicates that the function \eqref{eq: yamate} has local maximums and local minimums 
whose absolute values are exactly the same.

Tables \ref{table: E-optimal designs w = 1}--\ref{table: E-optimal designs w = 1 - x 1 + x} in Appendix \ref{sec: numerical example of dette algorithm} give some numerical examples.
We can confirm that E-optimal designs for weighted regression with weight functions \eqref{eq: dette weight function}
can be computed by using Theorem \ref{thm: dette E-optimal design for special weight function}.

\subsection{Approximate Tchebycheff Functions and An Approach to E-optimal Designs for General Weight Functions} 
\label{sec: approximate Tchebycheff system and approximate E-optimal design}

By Theorem \ref{thm: dette E-optimal design for special weight function}, 
if the weight function $w(x)$ of regression is given by \eqref{eq: dette weight function}, 
then the E-optimal designs can be computed. 
In this subsection, we discuss the weighted polynomial regression of general weight functions 
and approximate approach to it. 

In the previous sections, we discuss the exact Tchebycheff designs.
We have only a few examples that we can obtain the Tchebycheff designs 
for weighted polynomial regression by the same way as in Subsection \ref{sec: e-optimal design in dette weight function}. 
In this section we discuss a new relaxation of the Tchebycheff function $\kappa (x)$.
\begin{breakbox}
	\begin{defi}[Approximate Tchebycheff functions] \label{defi: approximate Tchebycheff function}
		Suppose that \\ $x \in [-1, 1]$.
		For a general weight function $w(x)$ of regression 
		such that if $-1 < x < 1$ then $w(x) > 0$, 
		the function $\kappa^\dagger (x)$ obtained by the following steps is called 
		an \keyword{approximate Tchebycheff function}.
		\begin{enumerate}[{\rm (a)}]
			\item Compute $v_m(x)$ by the Gram--Schmidt orthogonalization
			for the weight function $\eta (x)$ of orthogonal polynomials given by
			\begin{align}
				\eta (x) = \frac{w(x)}{\sqrt{1 - x^2}} , \label{eq: the relation between w and eta}
			\end{align}
			where $u_k (x) = x^{k - 1}, \ k = 1, 2, \dots , m$ are used in Gram--Schmidt orthogonalization.
			\item Obtain $\kappa^\dagger (x) = v_m(x) \sqrt{w(x)}$ .
		\end{enumerate}
		The Tchebycheff points $s_1^\dagger , s_2^\dagger , \dots , s_m^\dagger$ 
		of the approximate Tchebycheff function $\kappa^\dagger (x)$
		are defined as local maximum points and local minimum points.
	\end{defi}
\end{breakbox}
\begin{rema}[The difference of two ``weight functions'']
	Note that the weight function $w(x)$ of regression and 
	the weight function $\eta (x)$ of orthogonal polynomials are different.
	When the approximate Tchebycheff functions are considered,
	The relationship between $w(x)$ and $\eta (x)$ is denoted by \eqref{eq: the relation between w and eta}.
\end{rema}

\begin{figure}[!ht]
	\begin{center}
		\includegraphics[width=60mm, clip]{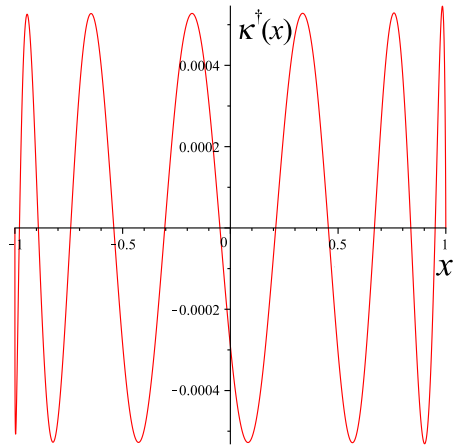}
		\caption{The approximate Tchebycheff function $\kappa^\dagger (x)$ of Definition \ref{defi: approximate Tchebycheff function}, 
		$m = 12, \ w(x) = (1 - x)^{1} \left( 3/2 + x \right)^{1/2}$} \label{fig: approximate Tpoly 3}
	\end{center}
\end{figure}

The exact Tchebycheff function $\kappa (x)$ has local maximums and local minimums
whose absolute values are exactly the same.
In contrast, the approximate Tchebycheff function $\kappa^\dagger (x)$ has local maximums and local minimums on $[-1, 1]$
whose absolute values are almost the same.
Figures~\ref{fig: approximate Tpoly 3} shows an example of the approximate Tchebycheff functions.
This figure indicates that an approximate Tchebycheff function $\kappa^\dagger (x)$ has 
local maximums and local minimums whose absolute values are almost the same but not necessary to be exactly the same.

By using Definition \ref{defi: approximate Tchebycheff function}, 
the proposed algorithm for calculating the approximation of E-optimal designs 
for weighted polynomial regression is described as follows.
\begin{breakbox}
	\begin{algo}[The algorithm for constructing the approximate Tchebycheff designs] \label{algo: algorithm for approximate Tchebycheff design}
		For a general weight function $w(x)$ of regression 
		such that if $-1 < x < 1$ then $w(x) > 0$, 
		the design $\mu^\dagger$ for weighted polynomial regression computed by the following steps
		is called the \keyword{approximate Tchebycheff design}.
		\begin{enumerate}[{\rm (a)}]
			\item Compute the approximate Tchebycheff function $\kappa^\dagger (x)$ by the Definition~\ref{defi: approximate Tchebycheff function}, 
			\item Compute Tchebycheff points $s_1^\dagger, s_2^\dagger, \dots , s_m^\dagger$ of the approximate Tchebycheff function $\kappa^\dagger (x)$, 
			\item Compute the design $\mu^\dagger$ given by
			\begin{align*}
				& \mu^\dagger = 
				\begin{pmatrix}
					s_1^\dagger & s_2^\dagger & \dots & s_m^\dagger \\
					\rho_1 & \rho_2 & \dots & \rho_m
				\end{pmatrix} , 
				\\
				& \left( \rho_1, \rho_2, \dots , \rho_m \right)^\top = \frac{F^{-1} \gamma}{\gamma^\top \gamma}
				, \quad F = \left( f_{i - 1}(s_j^\dagger ) \cdot (-1)^{j + 1} \right)_{i, j = 1, 2, \dots , m}
			\end{align*}
			where $\gamma = \left( \gamma_1, \gamma_2, \dots , \gamma_m \right)^\top$ denotes the vector of 
			the coefficients of the approximate Tchebycheff function $\kappa^\dagger (x)$.
		\end{enumerate}
	\end{algo}
\end{breakbox}

\begin{rema}[The approximate Tchebycheff functions on $\mathcal{X} = {[a, b]}$] 
	We can also consider the approximate Tchebycheff designs on $\mathcal{X} = [a, b]$ by the similar way.
	In order to calculate the approximate Tchebycheff functions for weighted polynomial regression 
	on $\mathcal{X} = [a, b]$ instead of $\mathcal{X} = [-1, 1]$, 
	the relationship between $w(x)$ and $\eta (x)$ 
	\begin{align*}
		\eta (x) = \frac{w(x)}{\sqrt{(x - a) (b - x)}}
	\end{align*}
	is used.
\end{rema}

In the next section, we show some results of numerical examples 
in order to verify that approximate Tchebycheff designs are close to E-optimal designs.

\section{Numerical Examples} \label{sec: numerical example}

In the previous section, we discuss Algorithm \ref{algo: algorithm for approximate Tchebycheff design},
the algorithm for constructing the approximate Tchebycheff designs for weighted polynomial regression.
In this section, we verify the accuracy of this algorithm with some numerical examples.

Tables \ref{table: approximate Tchebycheff designs w = 1-x12 2+x12}--\ref{table: approximate Tchebycheff designs w = 1 - x 1 + x} give numerical examples 
for Algorithm \ref{algo: algorithm for approximate Tchebycheff design}.
The first column contains the degrees $m$ of regression.
The second column contains the graphs of the approximate Tchebycheff functions $\kappa^\dagger (x)$.
The third column contains the experimental conditions $x_1, x_2, \dots , x_m$ and $\rho_1, \rho_2, \dots , \rho_m$
of the approximate Tchebycheff designs written as
\begin{align*}
	\mu = 
	\begin{pmatrix}
		x_1 & x_2 & \dots & x_m \\
		\rho_1 & \rho_2 & \dots & \rho_m
	\end{pmatrix}
	.
\end{align*}
The fourth column contains the minimum eigenvalues $\lamb_\mathrm{min} \left( M(\mu ) \right)$ of the Fisher information matrices.
The last column contains $1 - \mathrm{eff}_m^{\mathrm{E}} (\mu )$, 
where $\mathrm{eff}_m^{\mathrm{E}} (\mu )$ denotes \keyword{E-efficiency} of a design $\mu$ defined by 
\begin{align*}
	\mathrm{eff}_m^{\mathrm{E}} (\mu ) = \frac{\lamb_{\min} (M (\mu ))}{\underset{\mu}{\sup} \ \lamb_{\min} (M (\mu ))} .
\end{align*}
Note that the E-efficiency $\mathrm{eff}_m^{\mathrm{E}} (\mu ) \in [0, 1]$.
The E-efficiency $\mathrm{eff}_m^{\mathrm{E}}$ indicates that 
the larger the E-efficiency $\mathrm{eff}_m^{\mathrm{E}} (\mu )$ is, the better the design $\mu$ is in terms of the E-optimality criterion.
In order to compute the E-efficiency $\mathrm{eff}_m^{\mathrm{E}} (\mu )$, 
it is necessary to compute the E-optimal designs.
If $w(x) = (1 - x)^\alpha (1 + x)^\beta , \; \alpha , \beta \in \left\{ 0, 1 \right\}$, 
then the design is calculated by Theorem \ref{thm: dette E-optimal design for special weight function}.
Otherwise, the design calculated by the random optimization operated for a long time
is used instead of the E-optimal designs.
The computation is executed on the software Maple 15.

Tables \ref{table: approximate Tchebycheff designs w = 1-x12 2+x12}--\ref{table: approximate Tchebycheff designs w = 1 - x 1 + x} indicate that 
the E-efficiency $\mathrm{eff}_m^{\mathrm{E}} (\mu )$ is close to 1 regardless of the weight function $w(x)$.
This verifies the accuracy of Algorithm \ref{algo: algorithm for approximate Tchebycheff design}. 
Moreover, if the weight function $w(x)$ coincides \eqref{eq: dette weight function}, 
the E-efficiency $\mathrm{eff}_m^{\mathrm{E}} (\mu )$ is exactly equal to 1. 
This means that Algorithm \ref{algo: algorithm for approximate Tchebycheff design} computes the E-optimal designs exactly 
for weighted polynomial regression with the weight function $w(x) = (1 - x)^\alpha (1 + x)^\beta , \; \alpha , \beta \in \left\{ 0, 1 \right\}$.

\begin{table}[!htbp]
	\begin{center}
		\caption{Approximate Tchebycheff designs for weighted polynomial regression with $w(x) = (1 - x)^{1 / 2} (2 + x)^{1 / 2}$}
		\label{table: approximate Tchebycheff designs w = 1-x12 2+x12}
		\begin{tabular}{|c|c|l|c|c|} \hline 
			$m$ & graph of $\kappa^\dagger (x)$ & \multicolumn{1}{|c|}{appr. Tchebycheff design} & $\lamb_{\min} (M (\mu ))$ & $1 - \mathrm{eff}_m^{\mathrm{E}} (\mu )$ \\ \hline
			\multirow{5}{*}{$3$} & \multirow{5}{*}{\includegraphics[width=22mm, clip]{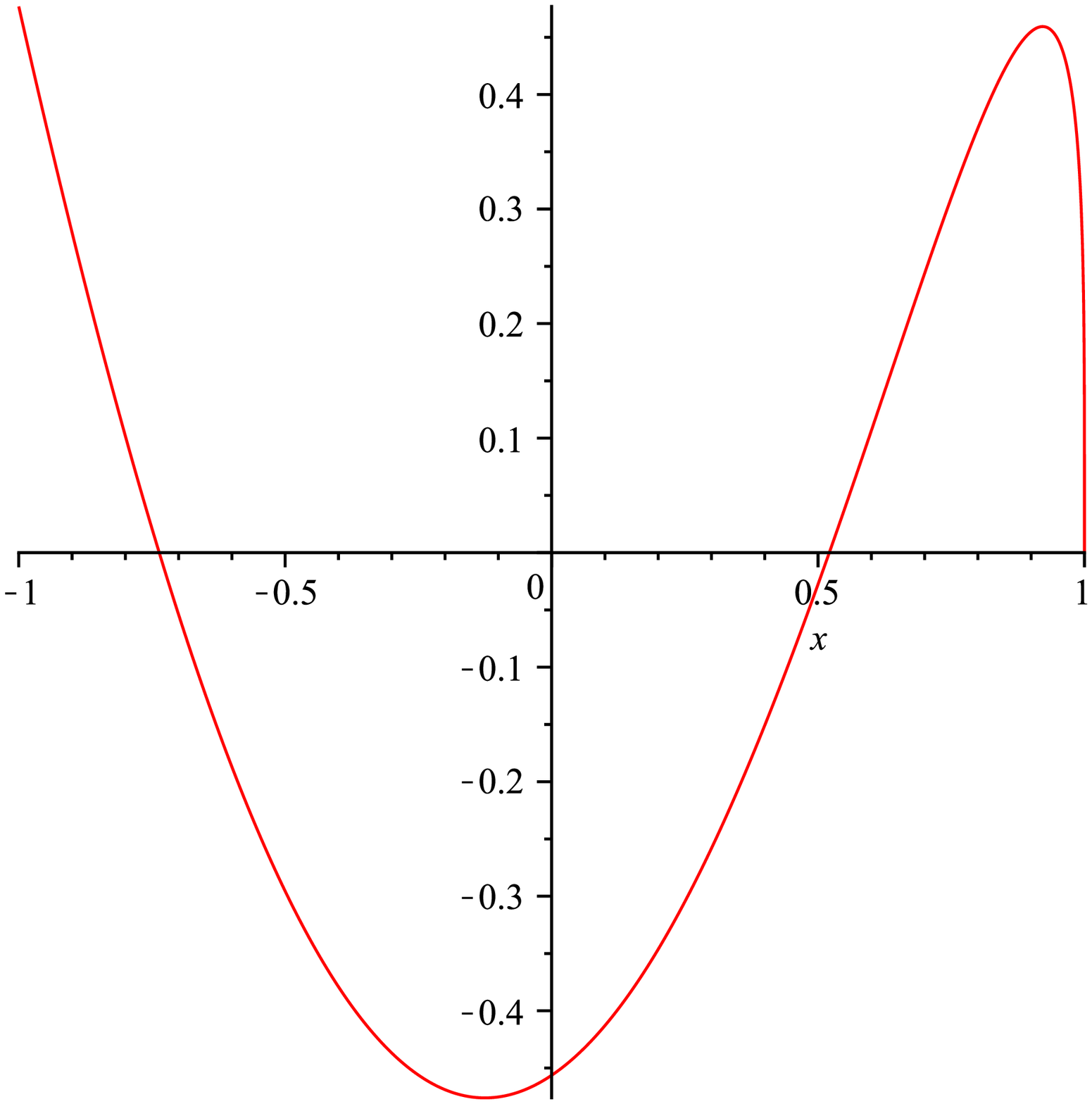}}
			      &                                    & \multirow{5}{*}{$7.693 \times 10^{-3}$} & \multirow{5}{*}{$8.720 \times 10^{-5}$} \\
			    & & $x_1 = -1.000, \ \rho_1 = 0.1721$ & & \\
			    & & $x_2 = -0.1252, \ \rho_2 = 0.4896$ & & \\
			    & & $x_3 = 0.9215, \ \rho_3 = 0.3383$  & & \\
			    & &                                    & & \\ \hline
			\multirow{10}{*}{$10$} & \multirow{10}{*}{\includegraphics[width=22mm, clip]{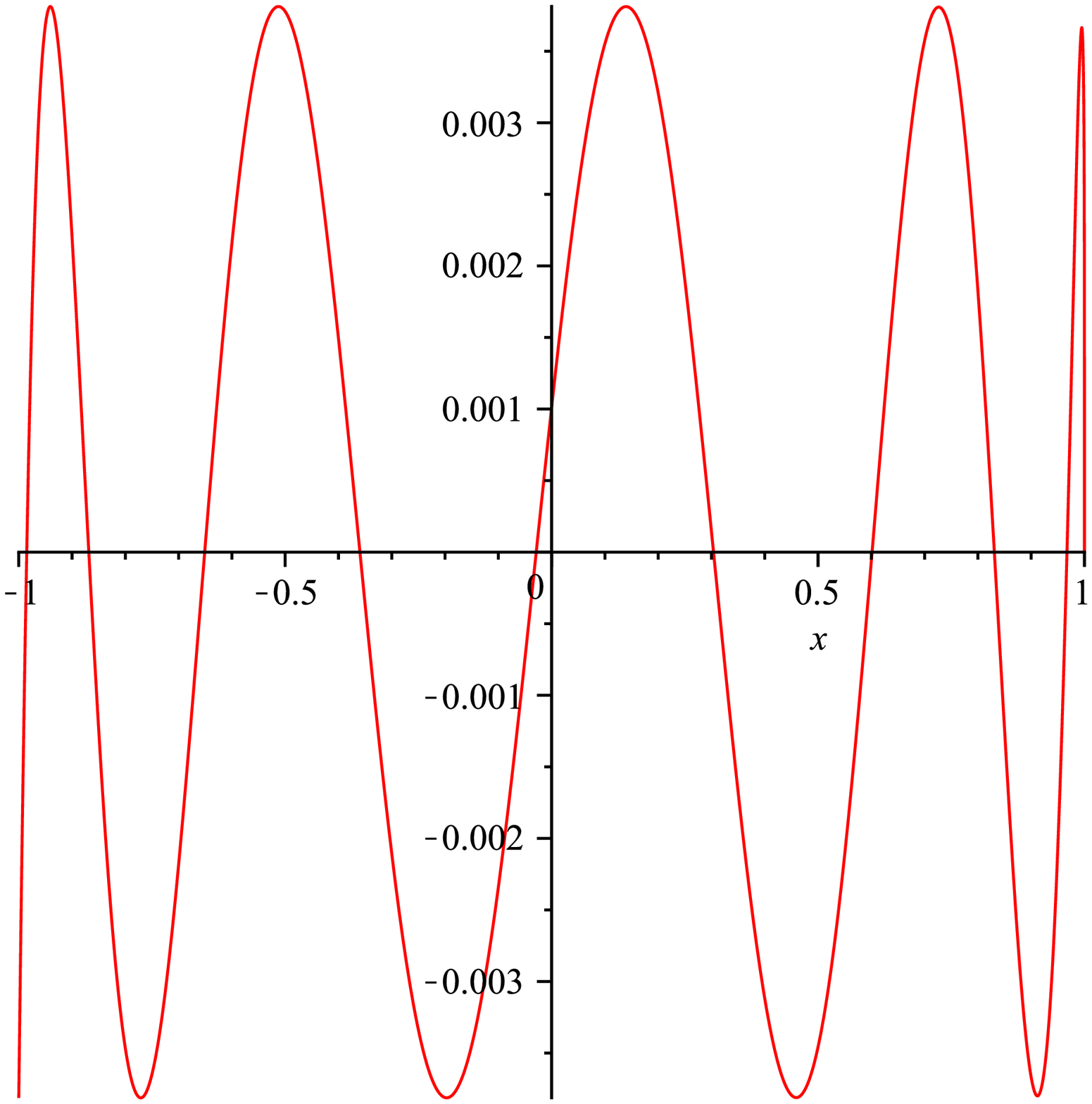}}
			      & $x_1 = -1.000, \ \rho_1 = 0.03909$ & \multirow{10}{*}{$1.714 \times 10^{-6}$} & \multirow{10}{*}{$3.334 \times 10^{-5}$} \\
			    & & $x_2 = -0.9407, \ \rho_2 = 0.08305$ & & \\
			    & & $x_3 = -0.7710, \ \rho_3 = 0.09785$ & & \\
			    & & $x_4 = -0.5126, \ \rho_4 = 0.1201$  & & \\
			    & & $x_5 = -0.1969, \ \rho_5 = 0.1395$ & & \\
			    & & $x_6 = 0.1396, \ \rho_6 = 0.1423$ & & \\
			    & & $x_7 = 0.4592, \ \rho_7 = 0.1261$ & & \\
			    & & $x_8 = 0.7269, \ \rho_8 = 0.1031$ & & \\
			    & & $x_9 = 0.9118, \ \rho_9 = 0.08509$ & & \\
			    & & $x_{10} = 0.9949, \ \rho_{10} = 0.06379$ & & \\ \hline
		\end{tabular}
		\caption{Approximate Tchebycheff designs for weighted polynomial regression with $w(x) = \ex^x$}
		\label{table: approximate Tchebycheff designs w = ex}
		\begin{tabular}{|c|c|l|c|c|} \hline 
			$m$ & graph of $\kappa^\dagger (x)$ & \multicolumn{1}{|c|}{appr. Tchebycheff design} & $\lamb_{\min} (M (\mu ))$ & $1 - \mathrm{eff}_m^{\mathrm{E}} (\mu )$ \\ \hline
			\multirow{5}{*}{$3$} & \multirow{5}{*}{\includegraphics[width=22mm, clip]{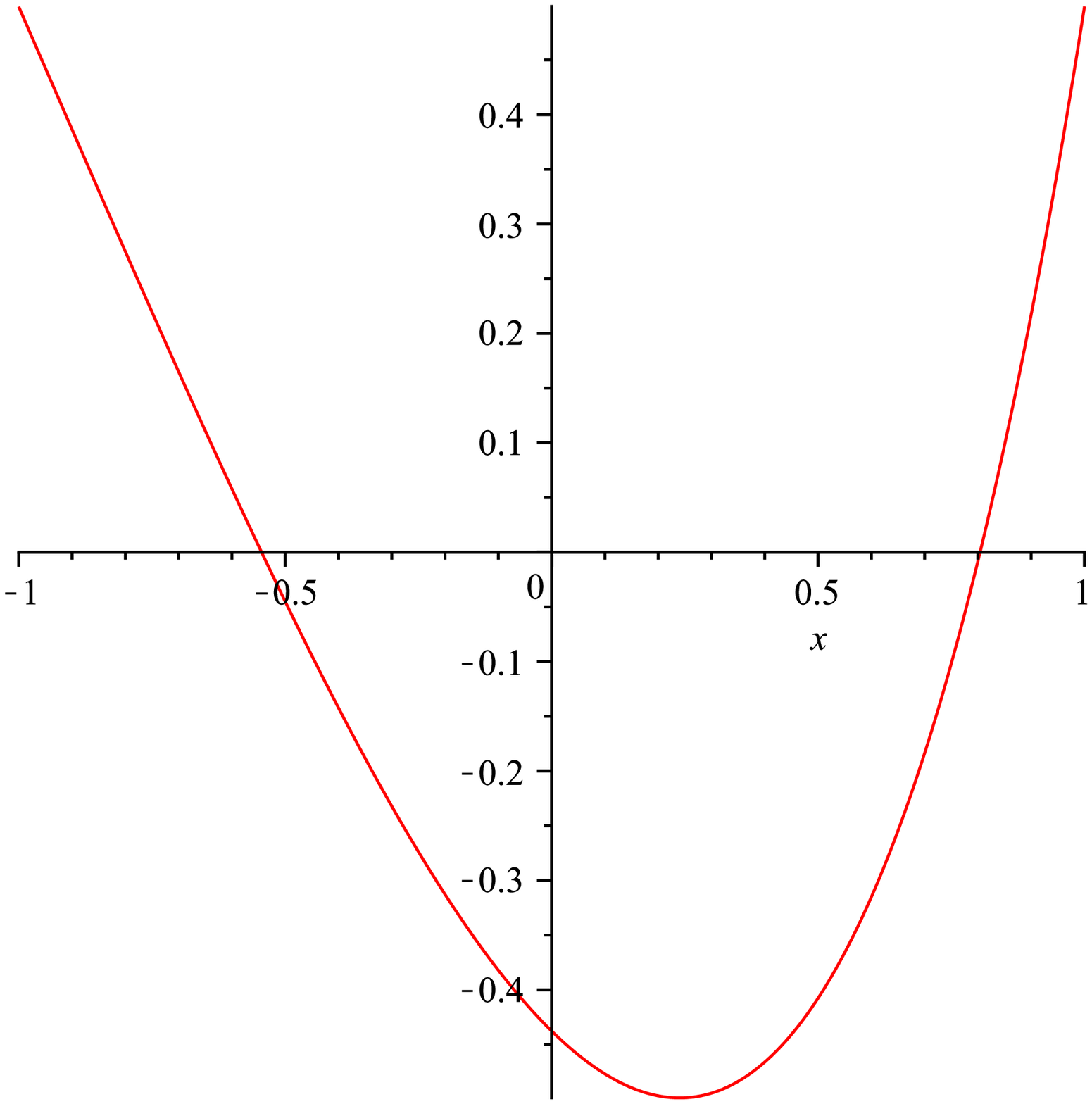}}
			      &                                    & \multirow{5}{*}{$1.976 \times 10^{-1}$} & \multirow{5}{*}{$4.082 \times 10^{-8}$} \\ 
			    & & $x_1 = -1.000, \ \rho_1 = 0.3204$ & & \\
			    & & $x_2 = 0.2405, \ \rho_2 = 0.5360$ & & \\
			    & & $x_3 = 1.000, \ \rho_3 = 0.1436$  & & \\
			    & &                                    & & \\ \hline
			\multirow{10}{*}{$10$} & \multirow{10}{*}{\includegraphics[width=22mm, clip]{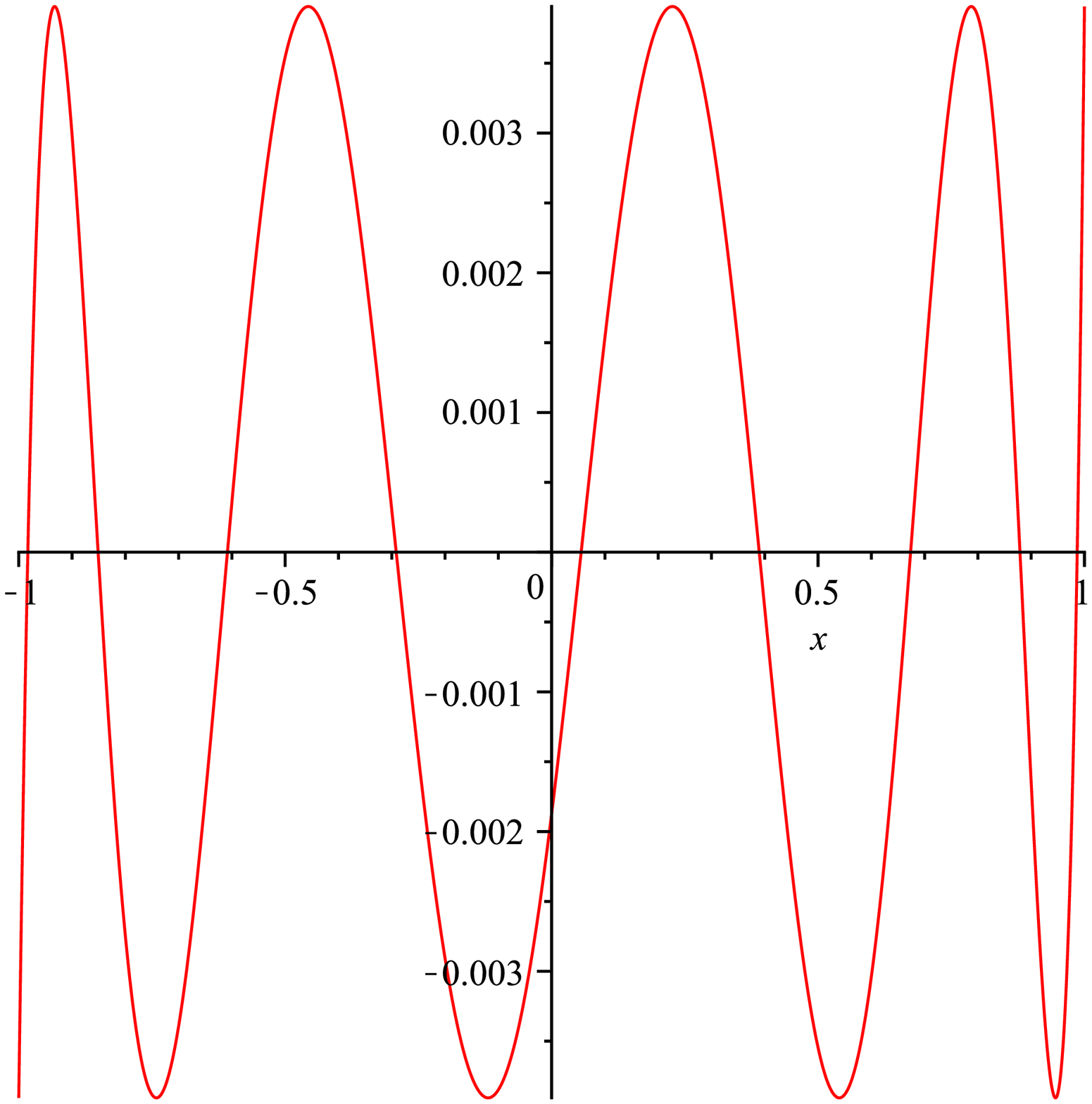}}
			      & $x_1 = -1.000, \ \rho_1 = 0.04351$ & \multirow{10}{*}{$1.660 \times 10^{-6}$} & \multirow{10}{*}{$2.998 \times 10^{-9}$} \\ 
			    & & $x_2 = -0.9326, \ \rho_2 = 0.09338$ & & \\
			    & & $x_3 = -0.7416, \ \rho_3 = 0.1119$ & & \\
			    & & $x_4 = -0.4566, \ \rho_4 = 0.1360$  & & \\
			    & & $x_5 = -0.1190, \ \rho_5 = 0.1494$ & & \\
			    & & $x_6 = 0.2267, \ \rho_6 = 0.1404$ & & \\
			    & & $x_7 = 0.5399, \ \rho_7 = 0.1164$ & & \\
			    & & $x_8 = 0.7876, \ \rho_8 = 0.09315$ & & \\
			    & & $x_9 = 0.9457, \ \rho_9 = 0.07880$ & & \\
			    & & $x_{10} = 1.000, \ \rho_{10} = 0.03710$ & & \\ \hline
		\end{tabular}
	\end{center}
\end{table}

\begin{table}[!htbp]
	\begin{center}
		\caption{Approximate Tchebycheff designs for weighted polynomial regression with $w(x) = 1$}
		\label{table: approximate Tchebycheff designs w = 1}
		\begin{tabular}{|c|c|l|c|c|} \hline 
			$m$ & graph of $\kappa^\dagger (x)$ & \multicolumn{1}{|c|}{appr. Tchebycheff design} & $\lamb_{\min} (M (\mu ))$ & $1 - \mathrm{eff}_m^{\mathrm{E}} (\mu )$ \\ \hline
			\multirow{5}{*}{$3$} & \multirow{5}{*}{\includegraphics[width=22mm, clip]{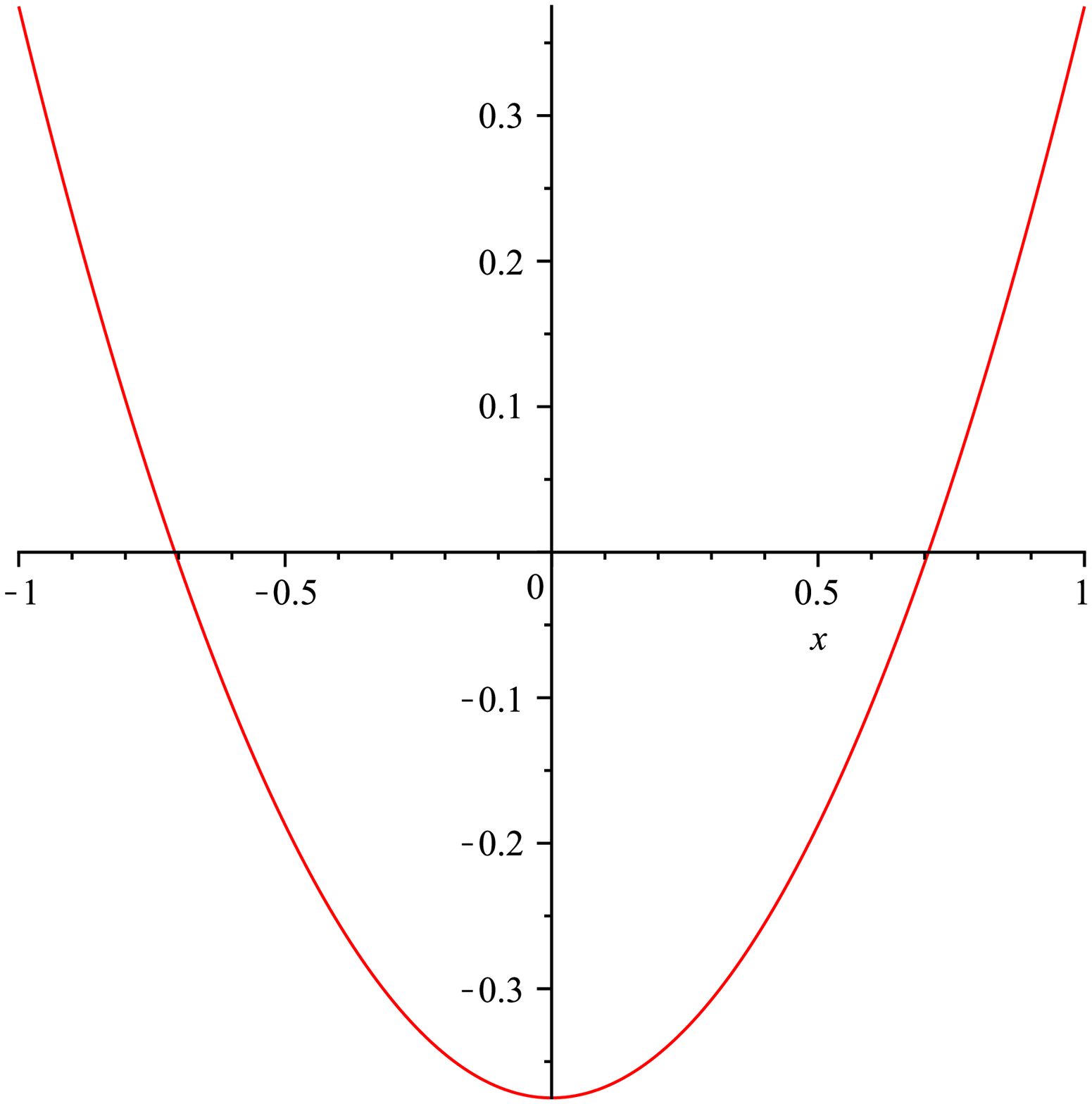}}
			      &                            & \multirow{5}{*}{$4.000 \times 10^{-1}$} & \multirow{5}{*}{$0.000$} \\
			    & & $x_1 = -1.000, \ \rho_1 = 0.2000$ &       & \\
			    & & $x_2 = 0.000, \ \rho_2 = 0.6000$  &       & \\ 
			    & & $x_3 = 1.000, \ \rho_3 = 0.2000$  &       & \\ 
			    & &                            &       & \\ \hline
			\multirow{10}{*}{$10$} & \multirow{10}{*}{\includegraphics[width=22mm, clip]{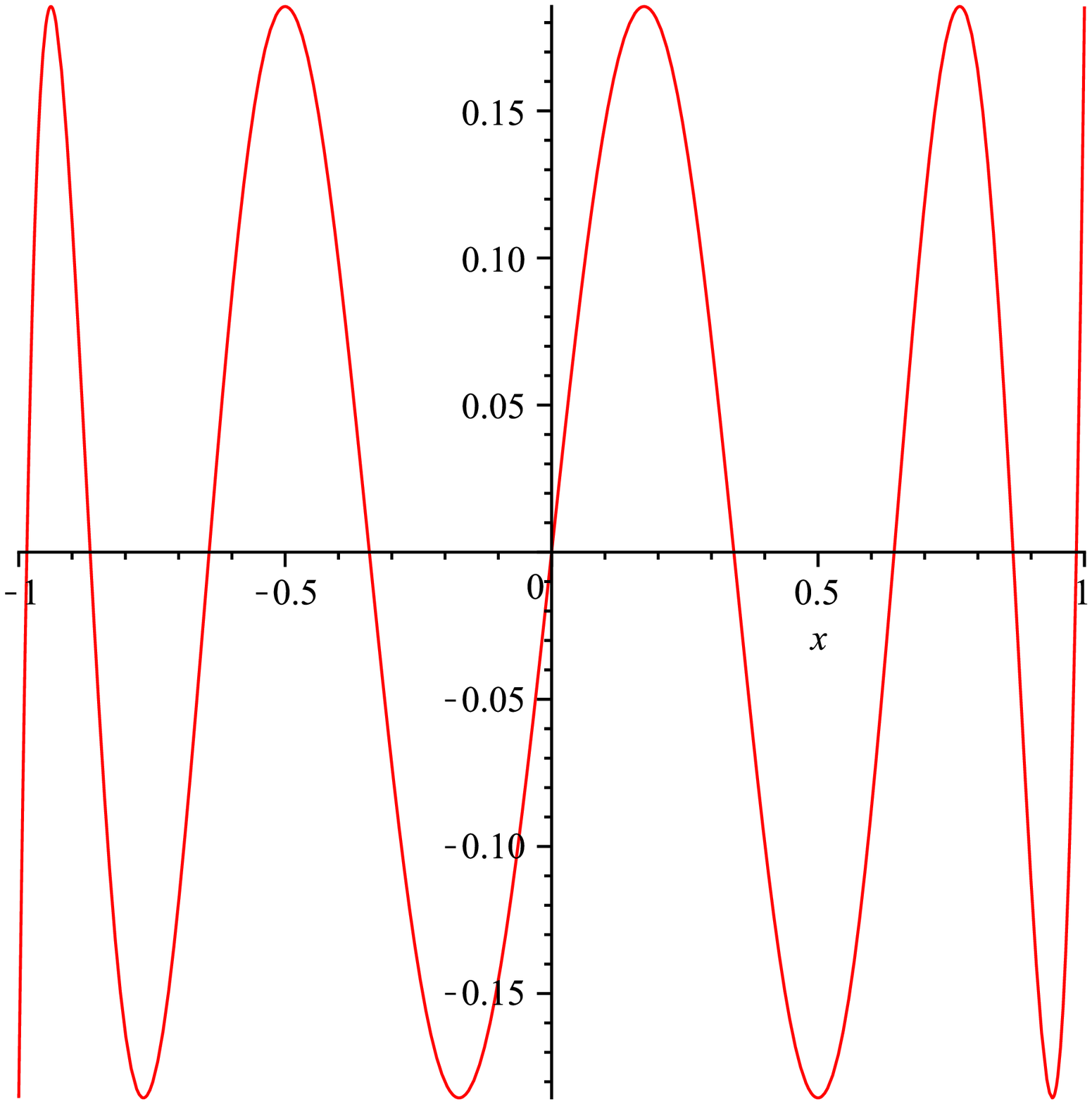}} 
			      & $x_1 = -1.000, \ \rho_1 = 0.04011$ & \multirow{10}{*}{$1.671 \times 10^{-6}$} & \multirow{10}{*}{$0.000$} \\
			    & & $x_2 = -0.9397, \ \rho_2 = 0.08563$ & & \\
			    & & $x_3 = -0.7660, \ \rho_3 = 0.1020$ & & \\
			    & & $x_4 = -0.5000, \ \rho_4 = 0.1263$ & & \\
			    & & $x_5 = -0.1736, \ \rho_5 = 0.1460$ & & \\
			    & & $x_6 = 0.1736, \ \rho_6 = 0.1460$ & & \\
			    & & $x_7 = 0.5000, \ \rho_7 = 0.1263$ & & \\
			    & & $x_8 = 0.7660, \ \rho_8 = 0.1020$ & & \\
			    & & $x_9 = 0.9397, \ \rho_9 = 0.08563$ & & \\
			    & & $x_{10} = 1.000, \ \rho_{10} = 0.04011$ & & \\ \hline
		\end{tabular}
	\end{center}
	\begin{center}
		\caption{Approximate Tchebycheff designs for weighted polynomial regression with $w(x) = 1 - x$}
		\label{table: approximate Tchebycheff designs w = 1 - x}
		\begin{tabular}{|c|c|l|c|c|} \hline 
			$m$ & graph of $\kappa^\dagger (x)$ & \multicolumn{1}{|c|}{appr. Tchebycheff design} & $\lamb_{\min} (M (\mu ))$ & $1 - \mathrm{eff}_m^{\mathrm{E}} (\mu )$ \\ \hline
			\multirow{5}{*}{$3$} & \multirow{5}{*}{\includegraphics[width=22mm, clip]{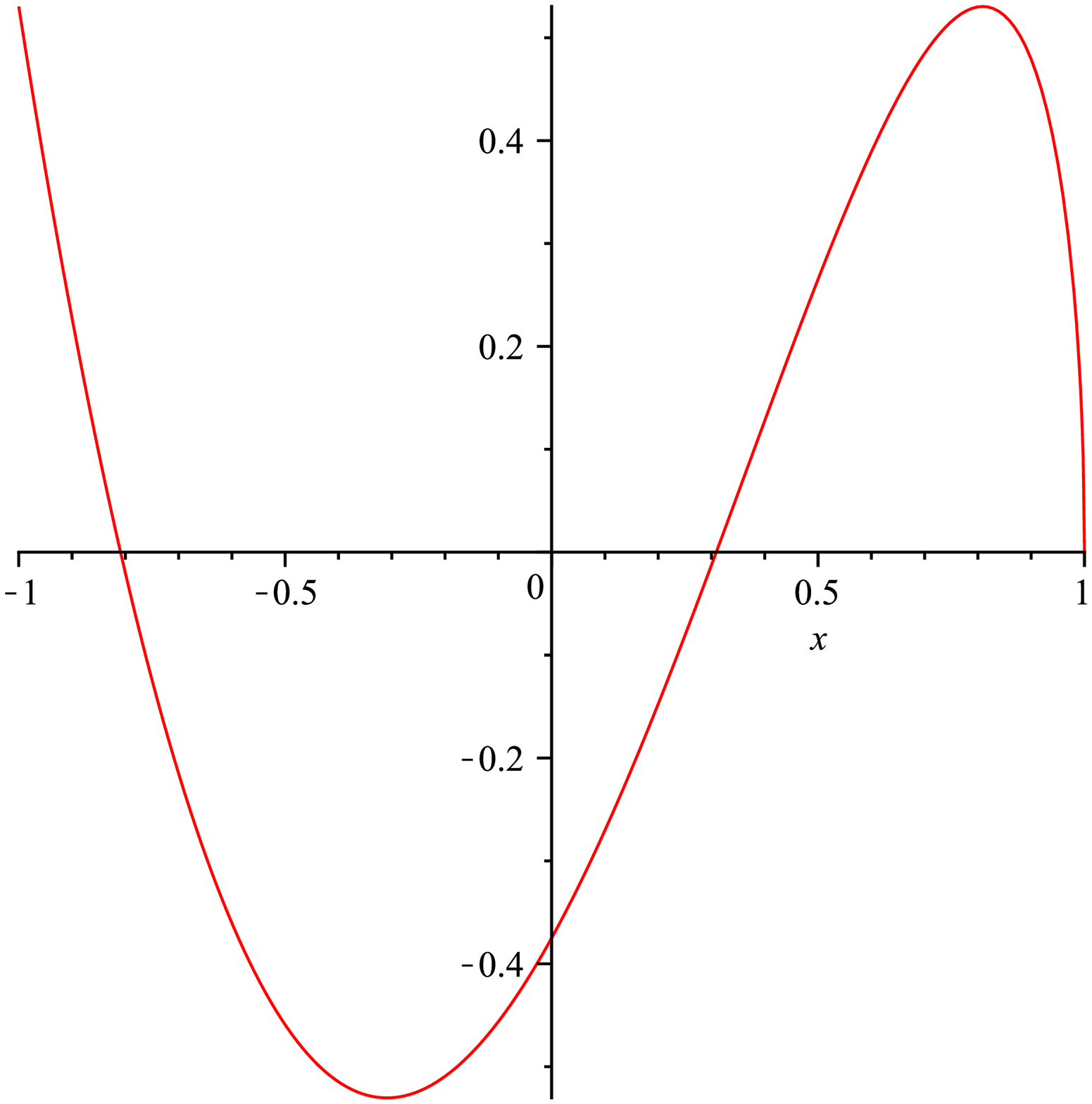}}
			      &                                    & \multirow{5}{*}{$9.524 \times 10^{-2}$} & \multirow{5}{*}{$0.000$} \\
			    & & $x_1 = -1.000, \ \rho_1 = 0.1238$      & & \\
			    & & $x_2 = -0.3090, \ \rho_2 = 0.3955$ & & \\
			    & & $x_3 = 0.8090, \ \rho_3 = 0.4807$  & & \\
			    & &                                    & & \\ \hline
			\multirow{10}{*}{$10$} & \multirow{10}{*}{\includegraphics[width=22mm, clip]{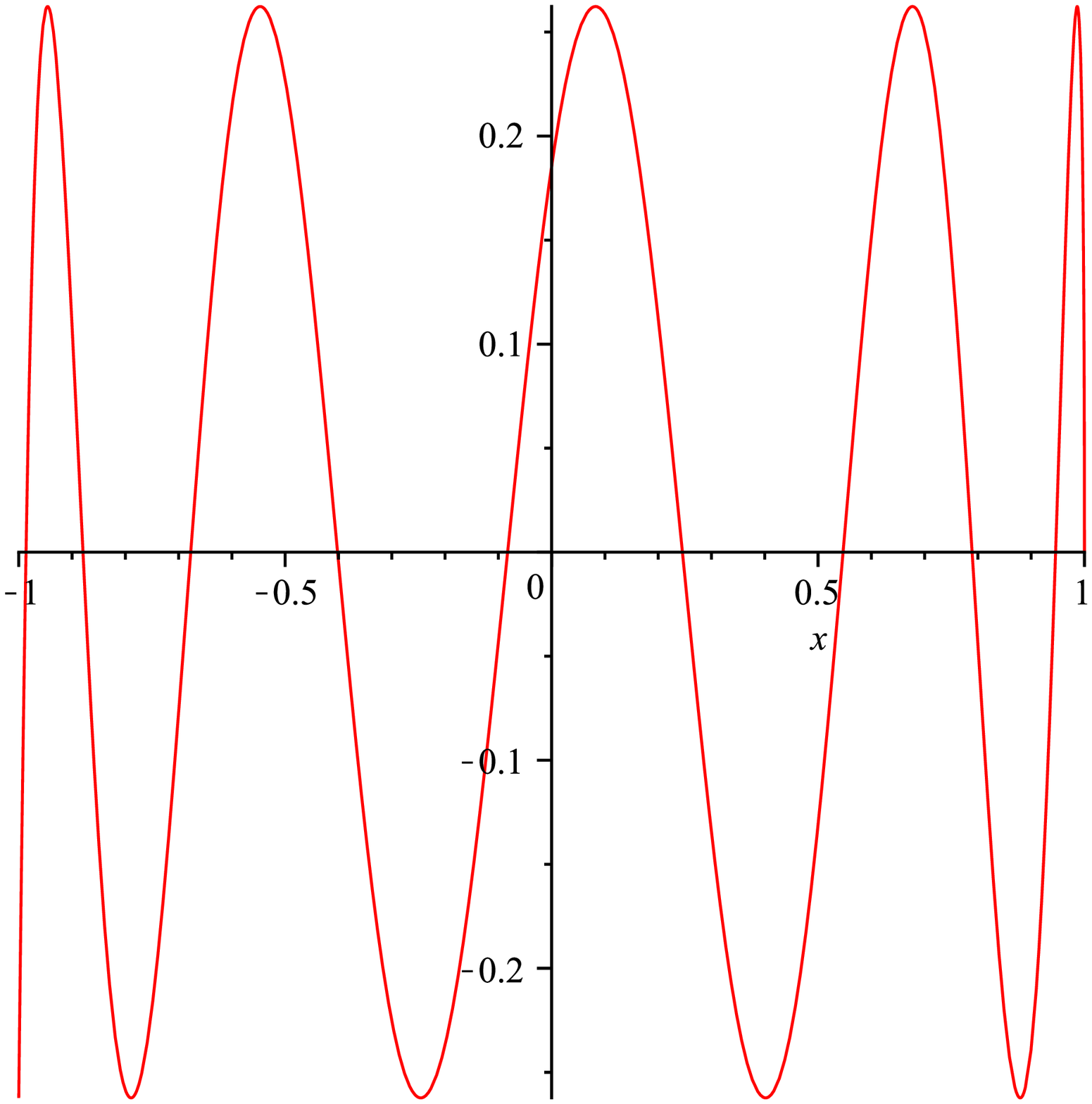}}
			      & $x_1 = -1.000, \ \rho_1 = 0.03642$     & \multirow{10}{*}{$9.463 \times 10^{-7}$} & \multirow{10}{*}{$0.000$} \\
			    & & $x_2 = -0.9458, \ \rho_2 = 0.07706$ & & \\
			    & & $x_3 = -0.7891, \ \rho_3 = 0.09006$ & & \\
			    & & $x_4 = -0.5469, \ \rho_4 = 0.1108$ & & \\
			    & & $x_5 = -0.2455, \ \rho_5 = 0.1321$ & & \\
			    & & $x_6 = 0.08258, \ \rho_6 = 0.1410$ & & \\
			    & & $x_7 = 0.4017, \ \rho_7 = 0.1311$ & & \\
			    & & $x_8 = 0.6773, \ \rho_8 = 0.1099$ & & \\
			    & & $x_9 = 0.8795, \ \rho_9 = 0.09082$ & & \\
			    & & $x_{10} = 0.9864, \ \rho_{10} = 0.08071$ & & \\ \hline
		\end{tabular}
	\end{center}
\end{table}

\begin{table}[!htbp]
	\begin{center}
		\caption{Approximate Tchebycheff designs for weighted polynomial regression with $w(x) = 1 + x$}
		\label{table: approximate Tchebycheff designs w = 1 + x}
		\begin{tabular}{|c|c|l|c|c|} \hline 
			$m$ & graph of $\kappa^\dagger (x)$ & \multicolumn{1}{|c|}{appr. Tchebycheff design} & $\lamb_{\min} (M (\mu ))$ & $1 - \mathrm{eff}_m^{\mathrm{E}} (\mu )$ \\ \hline
			\multirow{5}{*}{$3$} & \multirow{5}{*}{\includegraphics[width=22mm, clip]{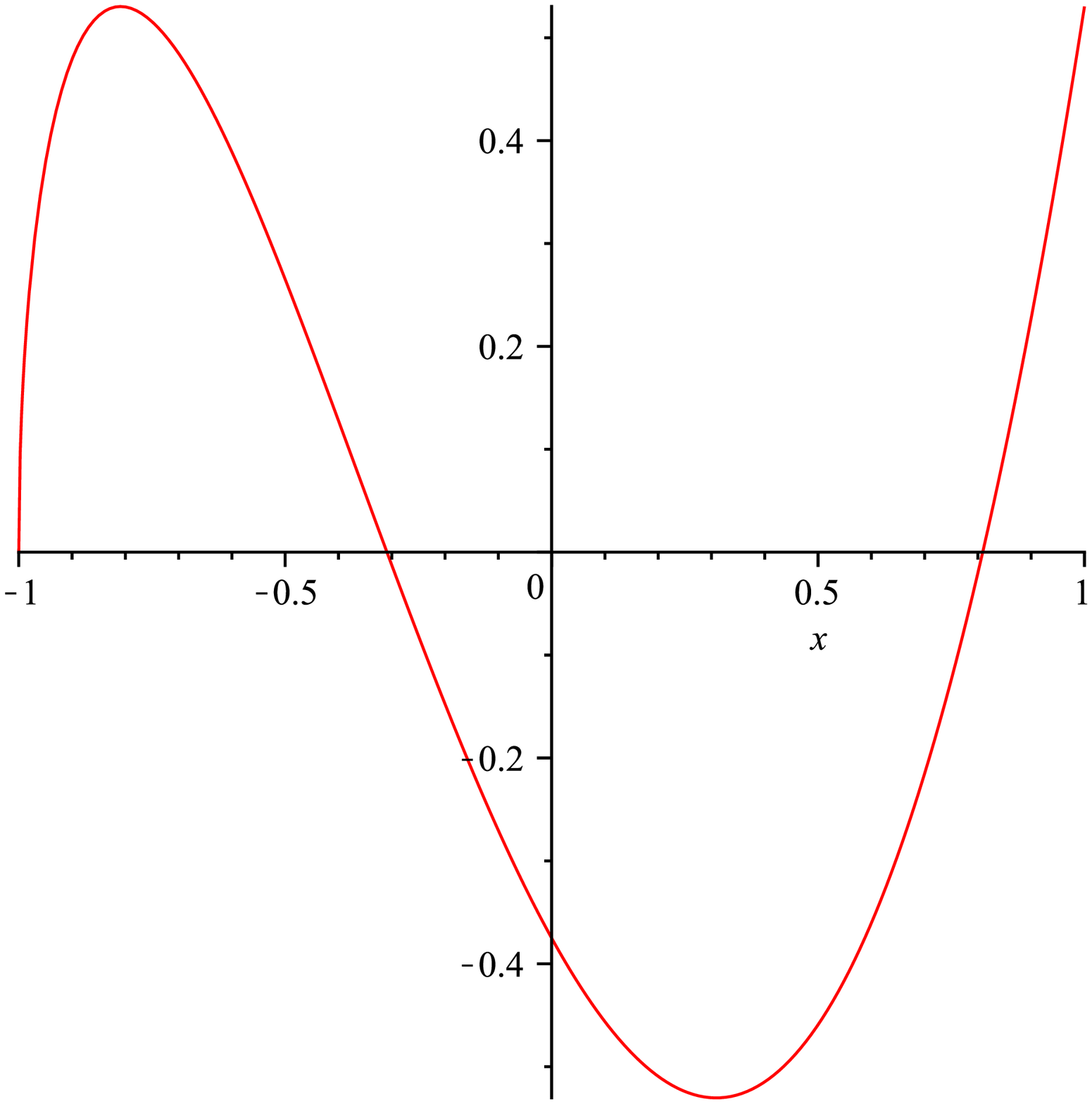}}
			      &                                    & \multirow{5}{*}{$9.524 \times 10^{-2}$} & \multirow{5}{*}{$0.000$} \\
			    & & $x_1 = -0.8090, \ \rho_1 = 0.4807$ & & \\
			    & & $x_2 = 0.3090, \ \rho_2 = 0.3955$ & & \\
			    & & $x_3 = 1.000, \ \rho_3 = 0.1238$  & & \\
			    & &                                    & & \\ \hline
			\multirow{10}{*}{$10$} & \multirow{10}{*}{\includegraphics[width=22mm, clip]{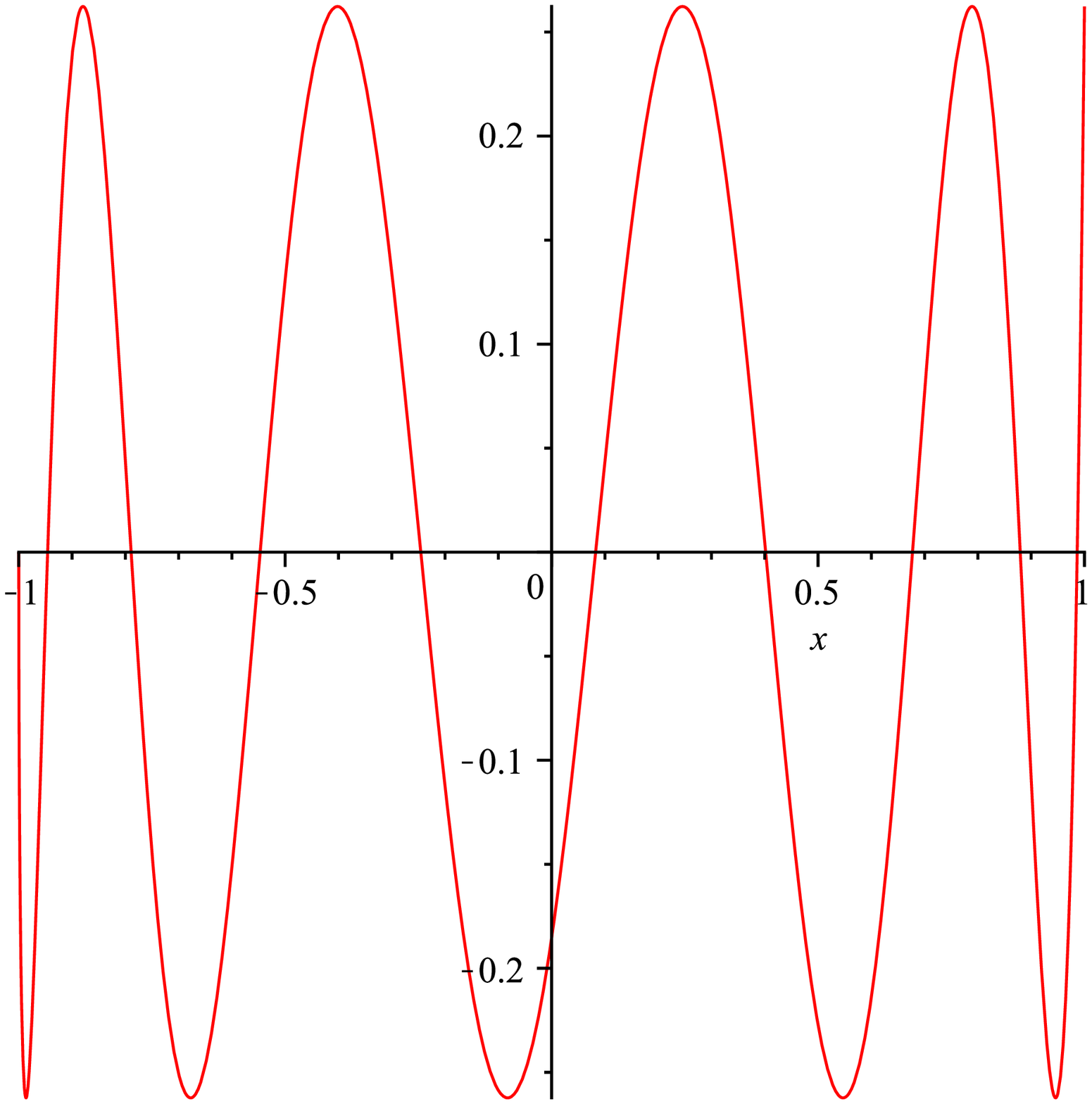}}
			      & $x_1 = -0.9864, \ \rho_1 = 0.08071$ & \multirow{10}{*}{$9.463 \times 10^{-7}$} & \multirow{10}{*}{$0.000$} \\
			    & & $x_2 = -0.8795, \ \rho_2 = 0.09082$ & & \\
			    & & $x_3 = -0.6773, \ \rho_3 = 0.01099$ & & \\
			    & & $x_4 = -0.4017, \ \rho_4 = 0.1311$  & & \\
			    & & $x_5 = -0.08258, \ \rho_5 = 0.1410$ & & \\
			    & & $x_6 = 0.2455, \ \rho_6 = 0.1321$ & & \\
			    & & $x_7 = 0.5469, \ \rho_7 = 0.1108$ & & \\
			    & & $x_8 = 0.7891, \ \rho_8 = 0.09006$ & & \\
			    & & $x_9 = 0.9458, \ \rho_9 = 0.07706$ & & \\
			    & & $x_{10} = 1.000, \ \rho_{10} = 0.03642$ & & \\ \hline
		\end{tabular}
	\end{center}
	\begin{center}
		\caption{Approximate Tchebycheff designs for weighted polynomial regression with $w(x) = (1 - x) (1 + x)$}
		\label{table: approximate Tchebycheff designs w = 1 - x 1 + x}
		\begin{tabular}{|c|c|l|c|c|} \hline 
			$m$ & graph of $\kappa^\dagger (x)$ & \multicolumn{1}{|c|}{appr. Tchebycheff design} & $\lamb_{\min} (M (\mu ))$ & $1 - \mathrm{eff}_m^{\mathrm{E}} (\mu )$ \\ \hline
			\multirow{5}{*}{$3$} & \multirow{5}{*}{\includegraphics[width=22mm, clip]{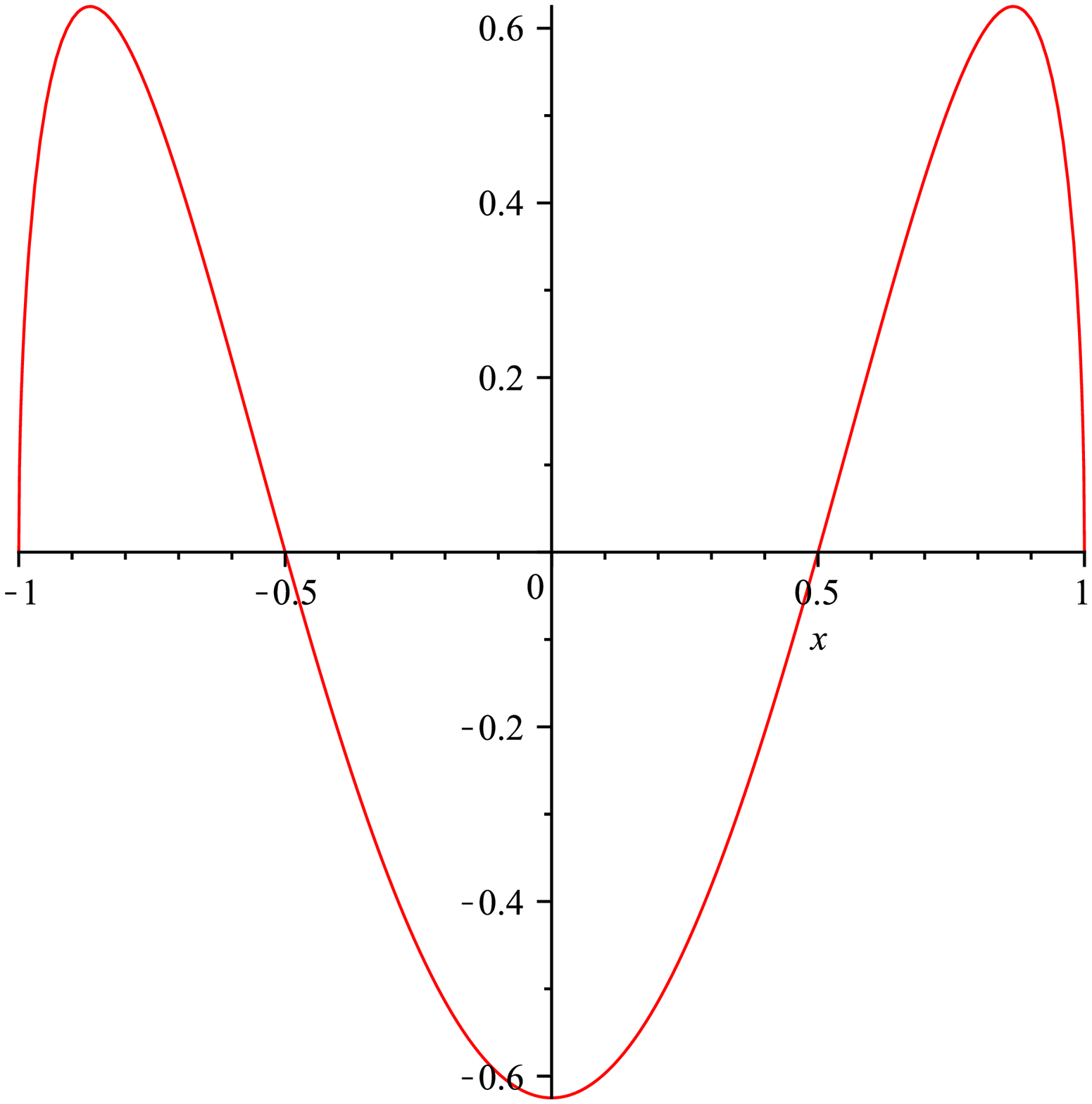}}
			      &                            & \multirow{5}{*}{$5.882 \times 10^{-2}$} & \multirow{5}{*}{$0.000$} \\
			    & & $x_1 = -0.8660, \ \rho_1 = 0.3137$ &       & \\
			    & & $x_2 = 0.000, \ \rho_2 = 0.3725$       &       & \\ 
			    & & $x_3 = 0.8660, \ \rho_3 = 0.3137$  &       & \\ 
			    & &                            &       & \\ \hline
			\multirow{10}{*}{$10$} & \multirow{10}{*}{\includegraphics[width=22mm, clip]{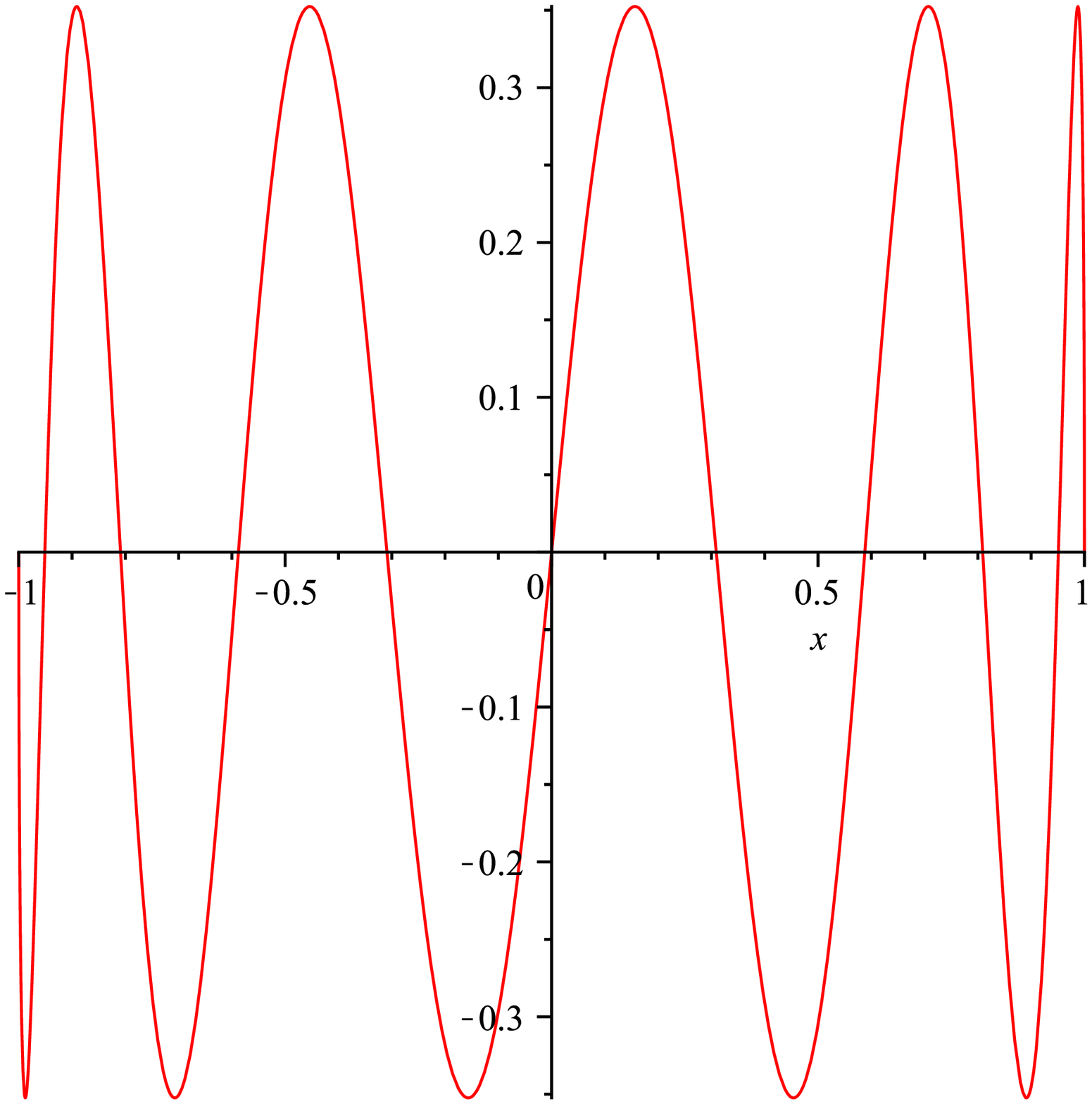}} 
			      & $x_1 = -0.9877, \ \rho_1 = 0.07329$ & \multirow{10}{*}{$5.593 \times 10^{-7}$} & \multirow{10}{*}{$0.000$} \\
			    & & $x_2 = -0.8910, \ \rho_2 = 0.08127$ & & \\
			    & & $x_3 = -0.7071, \ \rho_3 = 0.09702$ & & \\
			    & & $x_4 = -0.4540, \ \rho_4 = 0.1169$ & & \\
			    & & $x_5 = -0.1564, \ \rho_5 = 0.1315$ & & \\
			    & & $x_6 = 0.1564, \ \rho_6 = 0.1315$ & & \\
			    & & $x_7 = 0.4540, \ \rho_7 = 0.1169$ & & \\
			    & & $x_8 = 0.7071, \ \rho_8 = 0.09702$ & & \\
			    & & $x_9 = 0.8910, \ \rho_9 = 0.08127$ & & \\
			    & & $x_{10} = 0.9877, \ \rho_{10} = 0.07329$ & & \\ \hline
		\end{tabular}
	\end{center}
\end{table}

\clearpage

\section{Concluding Remarks} \label{sec: concluding}

In this thesis, we first indicate a new definition of approximate Tchebycheff functions.
By using this definition, we propose a new algorithm for
constructing the approximate Tchebycheff designs 
for weighted polynomial regression with general weight functions.
After that, we verify that the approximate Tchebycheff designs are close to E-optimal designs
by numerical examples.
Namely, the proposed algorithm enables us to calculate E-optimal designs approximately
for weighted polynomial regression with general weight functions.

As future works, it is necessary to discuss the definition of approximate Tchebycheff functions more strictly.
We must clarify how much gaps of the absolute values of local maximums and local minimums 
of approximate Tchebycheff functions are admitted.
Moreover, we must discuss how the proposed algorithm has a wide application.
We wish that the proposed algorithm is used to solve real problems of experimental designs.

%

\addcontentsline{toc}{section}{\quad \ \ References}


\appendix

\section{Appendix: Results of Numerical Examples}

\subsection{The Numerical Examples of Theorem \ref{thm: dette E-optimal design for special weight function}}
\label{sec: numerical example of dette algorithm}

Tables \ref{table: E-optimal designs w = 1}--\ref{table: E-optimal designs w = 1 - x 1 + x} are the results of the numerical examples for Theorem \ref{thm: dette E-optimal design for special weight function}.
The first column contains the degrees $m$ of regression.
The second column contains the graphs of Tchebycheff functions $\kappa (x)$.
The third column contains the experimental conditions $x_1, x_2, \dots , x_m$ and $\rho_1, \rho_2, \dots , \rho_m$
of the E-optimal designs written as
\begin{align*}
	\mu = 
	\begin{pmatrix}
		x_1 & x_2 & \dots & x_m \\
		\rho_1 & \rho_2 & \dots & \rho_m
	\end{pmatrix}
	.
\end{align*}
The last column contains the minimum eigenvalues $\lamb_\mathrm{min} \left( M(\mu ) \right)$ of the Fisher information matrix.
The computation is executed on the software Maple 15.

\begin{table}[!htbp]
	\begin{center}
		\caption{E-optimal designs for weighted polynomial regression with $w(x) = 1$}
		\label{table: E-optimal designs w = 1}
		\begin{tabular}{|c|c|l|c|} \hline 
			$m$ & graph of $\kappa (x)$ & \multicolumn{1}{|c|}{E-optimal design} & $\lamb_{\min} (M (\mu ))$ \\ \hline
			\multirow{5}{*}{$3$} & \multirow{5}{*}{\includegraphics[width=22mm, clip]{./jacobiP_Tpoly_00_n3.eps}}
			      &                            & \multirow{5}{*}{$4.000 \times 10^{-1}$} \\
			    & & $x_1 = -1.000, \ \rho_1 = 0.2000$ &       \\
			    & & $x_2 = 0.000, \ \rho_2 = 0.6000$  &       \\ 
			    & & $x_3 = 1.000, \ \rho_3 = 0.2000$  &       \\ 
			    & &                            &       \\ \hline
			\multirow{10}{*}{$10$} & \multirow{10}{*}{\includegraphics[width=22mm, clip]{./jacobiP_Tpoly_00_n10.eps}} 
			      & $x_1 = -1.000, \ \rho_1 = 0.04011$ & \multirow{10}{*}{$1.671 \times 10^{-6}$} \\
			    & & $x_2 = -0.9397, \ \rho_2 = 0.08563$ & \\
			    & & $x_3 = -0.7660, \ \rho_3 = 0.1020$ & \\
			    & & $x_4 = -0.5000, \ \rho_4 = 0.1263$ & \\
			    & & $x_5 = -0.1736, \ \rho_5 = 0.1460$ & \\
			    & & $x_6 = 0.1736, \ \rho_6 = 0.1460$ & \\
			    & & $x_7 = 0.5000, \ \rho_7 = 0.1263$ & \\
			    & & $x_8 = 0.7660, \ \rho_8 = 0.1020$ & \\
			    & & $x_9 = 0.9397, \ \rho_9 = 0.08563$ & \\
			    & & $x_{10} = 1.000, \ \rho_{10} = 0.04011$ & \\ \hline
		\end{tabular}
	\end{center}
	\begin{center}
		\caption{E-optimal designs for weighted polynomial regression with $w(x) = 1 - x$}
		\label{table: E-optimal designs w = 1 - x}
		\begin{tabular}{|c|c|l|c|} \hline 
			$m$ & graph of $\kappa (x)$ & \multicolumn{1}{|c|}{E-optimal design} & $\lamb_{\min} (M (\mu ))$ \\ \hline
			\multirow{5}{*}{$3$} & \multirow{5}{*}{\includegraphics[width=22mm, clip]{./jacobiP_Tpoly_10_n3.eps}}
			      &                                    & \multirow{5}{*}{$9.524 \times 10^{-2}$} \\
			    & & $x_1 = -1.000, \ \rho_1 = 0.1238$      & \\
			    & & $x_2 = -0.3090, \ \rho_2 = 0.3955$ & \\
			    & & $x_3 = 0.8090, \ \rho_3 = 0.4807$  & \\
			    & &                                    & \\ \hline
			\multirow{10}{*}{$10$} & \multirow{10}{*}{\includegraphics[width=22mm, clip]{./jacobiP_Tpoly_10_n10.eps}}
			      & $x_1 = -1.000, \ \rho_1 = 0.03642$     & \multirow{10}{*}{$9.463 \times 10^{-7}$} \\
			    & & $x_2 = -0.9458, \ \rho_2 = 0.07706$ & \\
			    & & $x_3 = -0.7891, \ \rho_3 = 0.09006$ & \\
			    & & $x_4 = -0.5469, \ \rho_4 = 0.1108$ & \\
			    & & $x_5 = -0.2455, \ \rho_5 = 0.1321$ & \\
			    & & $x_6 = 0.08258, \ \rho_6 = 0.1410$ & \\
			    & & $x_7 = 0.4017, \ \rho_7 = 0.1311$ & \\
			    & & $x_8 = 0.6773, \ \rho_8 = 0.1099$ & \\
			    & & $x_9 = 0.8795, \ \rho_9 = 0.09082$ & \\
			    & & $x_{10} = 0.9864, \ \rho_{10} = 0.08071$ & \\ \hline
		\end{tabular}
	\end{center}
\end{table}

\begin{table}[!htbp]
	\begin{center}
		\caption{E-optimal designs for weighted polynomial regression with $w(x) = 1 + x$}
		\label{table: E-optimal designs w = 1 + x}
		\begin{tabular}{|c|c|l|c|} \hline 
			$m$ & graph of $\kappa (x)$ & \multicolumn{1}{|c|}{E-optimal design} & $\lamb_{\min} (M (\mu ))$ \\ \hline
			\multirow{5}{*}{$3$} & \multirow{5}{*}{\includegraphics[width=22mm, clip]{./jacobiP_Tpoly_01_n3.eps}}
			      &                                    & \multirow{5}{*}{$9.524 \times 10^{-2}$} \\
			    & & $x_1 = -0.8090, \ \rho_1 = 0.4807$ & \\
			    & & $x_2 = 0.3090, \ \rho_2 = 0.3955$ & \\
			    & & $x_3 = 1.000, \ \rho_3 = 0.1238$  & \\
			    & &                                    & \\ \hline
			\multirow{10}{*}{$10$} & \multirow{10}{*}{\includegraphics[width=22mm, clip]{./jacobiP_Tpoly_01_n10.eps}}
			      & $x_1 = -0.9864, \ \rho_1 = 0.08071$ & \multirow{10}{*}{$9.463 \times 10^{-7}$} \\
			    & & $x_2 = -0.8795, \ \rho_2 = 0.09082$ & \\
			    & & $x_3 = -0.6773, \ \rho_3 = 0.01099$ & \\
			    & & $x_4 = -0.4017, \ \rho_4 = 0.1311$  & \\
			    & & $x_5 = -0.08258, \ \rho_5 = 0.1410$ & \\
			    & & $x_6 = 0.2455, \ \rho_6 = 0.1321$ & \\
			    & & $x_7 = 0.5469, \ \rho_7 = 0.1108$ & \\
			    & & $x_8 = 0.7891, \ \rho_8 = 0.09006$ & \\
			    & & $x_9 = 0.9458, \ \rho_9 = 0.07706$ & \\
			    & & $x_{10} = 1.000, \ \rho_{10} = 0.03642$ & \\ \hline
		\end{tabular}
	\end{center}
	\begin{center}
		\caption{E-optimal designs for weighted polynomial regression with $w(x) = (1 - x) (1 + x)$}
		\label{table: E-optimal designs w = 1 - x 1 + x}
		\begin{tabular}{|c|c|l|c|} \hline 
			$m$ & graph of $\kappa (x)$ & \multicolumn{1}{|c|}{E-optimal design} & $\lamb_{\min} (M (\mu ))$ \\ \hline
			\multirow{5}{*}{$3$} & \multirow{5}{*}{\includegraphics[width=22mm, clip]{./jacobiP_Tpoly_11_n3.eps}}
			      &                            & \multirow{5}{*}{$5.882 \times 10^{-2}$} \\
			    & & $x_1 = -0.8660, \ \rho_1 = 0.3137$ &       \\
			    & & $x_2 = 0.000, \ \rho_2 = 0.3725$       &       \\ 
			    & & $x_3 = 0.8660, \ \rho_3 = 0.3137$  &       \\ 
			    & &                            &       \\ \hline
			\multirow{10}{*}{$10$} & \multirow{10}{*}{\includegraphics[width=22mm, clip]{./jacobiP_Tpoly_11_n10.eps}} 
			      & $x_1 = -0.9877, \ \rho_1 = 0.07329$ & \multirow{10}{*}{$5.593 \times 10^{-7}$} \\
			    & & $x_2 = -0.8910, \ \rho_2 = 0.08127$ & \\
			    & & $x_3 = -0.7071, \ \rho_3 = 0.09702$ & \\
			    & & $x_4 = -0.4540, \ \rho_4 = 0.1169$ & \\
			    & & $x_5 = -0.1564, \ \rho_5 = 0.1315$ & \\
			    & & $x_6 = 0.1564, \ \rho_6 = 0.1315$ & \\
			    & & $x_7 = 0.4540, \ \rho_7 = 0.1169$ & \\
			    & & $x_8 = 0.7071, \ \rho_8 = 0.09702$ & \\
			    & & $x_9 = 0.8910, \ \rho_9 = 0.08127$ & \\
			    & & $x_{10} = 0.9877, \ \rho_{10} = 0.07329$ & \\ \hline
		\end{tabular}
	\end{center}
\end{table}

\end{document}